\documentstyle[12pt]{article}
\newcommand{\en}{\enspace}
\newcommand{\bi}{\bigskip}
\newcommand{\me}{\medskip}
\newcommand{\no}{\noindent}    
\newcommand{\non}{\nonumber}

\newcommand{\be}{\begin{equation}}                                              
\newcommand{\ee}{\end{equation}}                                                 
\newcommand{\bea}{\begin{eqnarray}}                                              
\newcommand{\eea}{\end{eqnarray}} 
\newcommand{\lra}{\longrightarrow}

\newcommand{\Z}{\makebox[0.06cm][l]{\sf Z}{\sf Z}}

\newcommand{\R}{{\rm I}\!{\rm R}}
\newcommand{\C}{{\mbox{\rm $\scriptscriptstyle ^\mid$\hspace{-0.40em}C}}} 

\renewcommand{\epsilon}{\varepsilon}

\title{A Classification Approach for Open Manifolds}
\author{J\"urgen Eichhorn, Greifswald}
\date{\today}

\newtheorem{satz}{Satz}[section]
\newtheorem{proposition}[satz]{Proposition}
\newtheorem{lemma}[satz]{Lemma}
\newtheorem{corollary}[satz]{Corollary}
\newtheorem{theorem}[satz]{Theorem}

\newcommand{\dil}{\mbox{dil }}
\newcommand{\diam}{\mbox{diam }}
\newcommand{\dist}{\mbox{dist }}
\newcommand{\diff}{\mbox{Diff }}
\newcommand{\intt}{\mbox{int }}
\newcommand{\all}{\mbox{all }}
\newcommand{\vol}{\mbox{vol }}

\newcommand{\comp}[1]{\mbox{\rm comp$_{#1}$\,}}
\newcommand{\sk}[2]{\mbox{$\langle#1,#2\rangle$}}

\begin{document}

\sloppy

\maketitle

\section{Introduction}

For closed manifolds there exists an effective highly
elaborated classification approach the main steps of which 
are the Thom--Pontrjagin construction, bordism theory, surgery,
Wall groups, the exact sequence of 
Browder--Novikov--Sullivan--Wall. All this can be expressed in an 
algebraic language of f. g. $Z \pi$--, unitary $Z \pi$--, 
$Q \pi$--modules and their $K$--theory etc. Moreover, 
there exist many number valued invariants like classical
characteristic numbers, signature, higher signature, analytic
and Reidemeister torsion, $\eta$--invariants. 

For open manifolds, absolutely nothing of this remains, 
at least at the first glance. We have the following simple

\begin{proposition}
Let ${\cal M}^n$ be the set of all smooth oriented manifolds
and $V$ a vector space or abelian group. There does not exist
a nontrivial map $c: {\cal M} \lra V$ such that

1. $M^n \cong M'$ orientation preserving diffeomorphic 
implies $c(M)=c(M')$ and

2. $c(M \# M')=c(M)+c(M')$.
\end{proposition}

\me
\no
{\bf Proof.}
Assume at first $M^n \not\cong \Sigma^n$, fix two points at $M^n$,
then $M_\infty = M_1 \# M_2 \# \dots$, $M_i (M,i) \cong M$
has a well defined meaning. We can write 
$M_\infty = M_1 \# M_{\infty, 2}$, 
$M_{\infty, 2} = M_2 \# M_3 \# \dots$
and get $c(M_\infty)=c(M)+c(M_{\infty,2})=c(M)+c(M_\infty)$,
$c(M)=0$. 

Assume $M^n=\Sigma^n$, and ord $\Sigma^n=k>1$ which yields
\[ c(\Sigma^n \# \dots \# \Sigma^n) = k \cdot c(\Sigma^n) =
   c(S^n), \quad c(\Sigma^n) = \frac{1}{k} c(S^n) , \]
\[ c(\Sigma^n) = c (\Sigma^n \# S^n) = (1+\frac{1}{k}) c(S^n).
   \quad c(S^n) = 2 c (S^n), \quad c(S^n) = 0 , \]
\[ c(\Sigma^n) = 0 . \]
\hfill $\Box$

The only number valued invariant defined for all connected
manifolds $M^n$ and known to the author is the dimension $n$.
If one characterizes orientabality / nonorientabality by
$\pm 1$ then there are two number--valued invariants. That
is all.

Denote by ${\cal M}^n([cl])$ the set of all diffeomorphism
classes of closed $n$--manifolds. Then we have

\begin{proposition}
$\# {\cal M}^n([cl]) = $ alef zero.
\end{proposition}

\me
\no
{\bf Proof.}
According to Cheeger, there are only finitely many
diffeomorphism types for $(M^n,g)$ with
$\diam (M^n,g) \le D$, $r_{inj} (M^n,g) \le i$, 
$|$ sectional curvature $(M^n,g)| \le K$, where
$r_{inj} (M^n,g)$ denotes the injectivity radius.
Setting $D_\nu = K_\nu = i_\nu = \nu$ and considering
$\nu \lra \infty$, we count all diffeomorphism types of
closed Riemannian $n$--manifolds, in particular all
diffeomorphism types of closed manifolds.
\hfill $\Box$

On the other hand for open manifolds holds

\begin{proposition}
The cardinality of ${\cal M}([open])$ is at least that
of the continuum, $n \ge 2$.
\end{proposition}

\me
\no
{\bf Proof.}
Assume $n \ge 3$, $n$ odd, let $2 = p_1 < p_2 < \dots$
be the increasing sequence of prime numbers and let
$L^n(p_\nu) = S^n / Z / p_\nu$ be the corresponding lens space.
Consider 
$M^n := d_1 \cdot L(p_1) \# d_2 \cdot L(p_2) \# \dots $,
$d_\nu=0,1$. Then any $0,1$--sequence $(d_1, d_2, \dots)$
defines a manifold and different sequences define non
diffeomorphic manifolds. If $n \ge 4$ is even multiply
with $S^1$. For $n=2$ the assertion follows from the
classification theorem in [20].
\hfill $\Box$

There are simple methods to construct only from one closed
manifold $M^n \neq \Sigma^n$ infinitely many nondiffeomorphic
manifolds. This, proposition 1.3 and other considerations
support the naive imagination, that ''measure of 
${\cal M}([open])$ : measure ${\cal M}([cl]) = \infty : 0$''.
We understand this as an additional hint how difficult
would be any classification of open manifolds.

A certain requirement of what should be the goal of such 
a classification comes from global analysis. The main task
of global analysis is the solution of linear and nonlinear
differential equations in dependence of the underlying
geometry and topology. Examples for open manifolds are
discussed and solved in [8], [9], [10]. To each open
manifold there are attached certain funcional spaces and 
admitted maps which enter into the classification should
induce maps between the functional spaces, i. e.
$f: M \lra M'$ induces a map $(f. sp.)(M') \lra (f. sp.)(M)$.
A very simple example shows that this is a reasonable
requirement. Consider the diffeomorphism
$f=\left( tg \left( \frac{\pi}{2} \cdot \right) \right)^{-1}: 
 ]0, \infty[ \lra ]0,1[$. Then $1 \in L_2(]0,1[)$
but $1=f^* 1 \notin L_2(]0, \infty[)$. Our first
conclusion is that one should classify pairs $(M,g)$ and
maps should be adapted to the Riemannian metrics under
consideration. A pure differential topological classification
will be to difficult, not handable and as we shortly indicated,
will be less important for applications. Manifolds which
appear in applications are endowed with a Riemannian metric.
Accepting this, the corresponding algebraic topology should
be analytic and simplicial $L_p$--(co)homology, bounded
(co)homology and others. Consider 
$(\R^n, g_{standard})=(\R^n, dr^2+r^2d\sigma^2_{S^{n-1}})$
and a canonical uniform triangulation $K_{\R^n}$ of
$\R^n$. Then 
$\overline{H}^{*,2}(\R^n) = H^{*,2}(\R^n) = H^{*,2}(K_{\R^n})
 = \overline{H}^{*,2}(K_{\R^n}) = 0$. 
But if we endow $\R^n$ with the hyperbolic metric
$g_H = dr^2+ (\sinh r)^2 d \sigma^2_{S^{n-1}}$
then for $n=2k$
$\dim \overline{H}^{k,2}(\R^n,g_H) = \dim H^{k,2}(\R^n,g_H) =
\infty$ and for $n=2k+1$ 
$\dim \overline{H}^{k,2}(\R^n,g_H) = 0$,
$\dim H^{k,2}(\R^n,g_H) = \dim H^{k+1,2}(\R^n,g_H) = \infty$.
Nevertheless $(\R^n, g_{st})$ and $(\R^n, g_H)$ are
canonical isomorphic by the best possible map $id_{\R^n}$.
We come to our preliminary conclusions,

1. admitted maps $(M,g) \lra (M',g')$ must be strongly
adapted to $g$, $g'$,

2. the classification approach will from a certain step on
be connected with spectral properties.

Our approach can be characterized as follows. We first 
decompose the set of all $(M,g)$'s into the set of 
components of a certain uniform structure and then classify
(up to a certain amount) the manifolds in the component under 
consideration, i. e. we have to define two classes of 
invariants, one for the components and one for the manifolds
inside a component. Moreover, we consider several uniform
structures which become finer and finer. This implies that
the (arc) components become smaller and smaller. Since the
whole approach is sufficiently extensive, we can present here
only the main steps and sketch some proofs.

The paper is organized as follows. In section 2, we define those
uniform structures which consider $(M^n,g)$ only as proper
metric spaces. The uniform structures in section 4 take into
account the smooth and Riemannian structure of $(M^n,g)$.
Section 5 is devoted to bordism theories, adapted in a certain
sense to uniform structures under consideration. In the
concluding section 6 we define and discuss classes of invariants
in the two senses sketched above. Details and complete proofs
are contained in [8], [9], [10]. All three papers will be
presented for publication.

This presentation is part of a talk given at the conference
in honor of V. A. Rohlin, St. Petersburg, August 1999.

\section{Uniform structures of proper metric spaces}

Let $Z=(Z,d_Z)$ be a metric spaces, $X, Y \subset Z$ subsets,
$\epsilon > 0$, define $U_\epsilon (X) = \{ z \in Z | 
\dist (z, X) < \epsilon \}$, analogously $U_\epsilon(Y)$.
Then the Hausdorff distance $d_H(X,Y)= d^Z_H (X,Y)$ is
defined as
\[ d^Z_H(X,Y) = \inf \{ \epsilon | X \subset U_\epsilon(Y),
   Y \subset U_\epsilon(X) \} . \]
If there is no such $\epsilon > 0$ then we set 
$d^Z_H (X,Y) = \infty$, $d^Z_H$ is an almost metric on the
set of all closed subsets, i. e. it has values in $[0,\infty]$
but satisfies all other conditions of a metric. If $Z$ ist
compact then $d^Z_H$ is a metric on the set of all closed 
subsets. A metric space $(X,d)$ is called proper if the closed
balls $\overline{B_\epsilon(x)}$ are compact for all
$x \in X$, $\epsilon >0$. This implies that $X$ is separable,
complete and locally compact. Any complete Riemannian manifold
$(M^n,g)$ is a proper metric space. In the sequel we restrict
to proper metric spaces. Let $(X,d_X)$, $(Y,d_Y)$ be metric
spaces, $X \sqcup Y$ their disjoint union. A metric on 
$X \sqcup Y$ is called admissible if $d$ restricts to $d_X$
and $d_Y$, respectively. The Gromov--Hausdorff distance
$d_{GH}(X,Y)$ is defined as
\[ d_{GH}(X,Y) = \inf \{ d^{X \sqcup Y}_H (X,Y) | d
   \mbox{ admissible on } X \sqcup Y \} . \]
Note that the Gromov--Hausdorff distance can be infinite. 
Gromov defined $d_{GH}$ originally as
\bea
   d_{GH}(X,Y) &=& \inf \{ d^Z_G (i(X),j(Y)) | i : X \lra Z, 
   j : Y \lra Z \mbox{ isometric embeddings} \non \\
   && \mbox{ into a metric space } Z \} . 
\eea 

\begin{lemma}
If $X$ and $Y$ are compact metric spaces and $d_{GH}(X,Y) = 0$
then $X$ and $Y$ are isometric.
\end{lemma}

\me
\no
{\bf Proof.}
This follows from the definition and an Arzela--Ascoli
argument.
\hfill $\Box$

Denote by ${\cal M}$ the set of all isometry classes
$[X]$ of proper metric spaces $X$ and 
${\cal M}_{GH} = {\cal M} / \sim$, where 
$[X] \sim [Y]$ if $d_{GH}([X],[Y])=0$.

\begin{proposition}
$d_{GH}$ defines an almost metric on ${\cal M}_{GH}$.
\end{proposition}
\hfill $\Box$

We denote in the sequel $X=[X]$ if there does not arise any
confusion.

Now we define the uniform structure. Let $\delta > 0$ and
set 
\[ V_\delta = \{ (X,Y) \in {\cal M}^2_{GH} | 
   d_{GH} (X,Y) < \epsilon \} . \]

\begin{lemma}
${\cal L} = \{ V_\delta \}_{\delta > 0}$ is a basis for a
metrizable uniform structure ${\cal U}_{GH} ({\cal M}_{GH})$.
\end{lemma}

\me
\no
{\bf Proof.}
${\cal L}$ is locally defined by a metric. Hence it satisfies all
desired conditions.
\hfill $\Box$

Let $\overline{\cal M}_{GH}$ be the completion of ${\cal M}_{GH}$ 
with respect to ${\cal U}_{GH}$ and denote the metric in 
$\overline{\cal M}_{GH}$ by $\overline{d}_{GH}$.

\begin{lemma}
${\cal M}_{GH} = \overline{\cal M}_{GH}$ as sets and $d_{GH}$ and
$\overline{d}_{GH}$ are locally equivalent.
\end{lemma}

\begin{proposition}
${\cal M}_{GH} = \overline{\cal M}_{GH}$ is locally arcwise
connected.
\end{proposition}

\me
\no
{\bf Proof.}
We refer to [8] for the proof which is quite elementary but
rather lengthy.
\hfill $\Box$

\begin{corollary}
In ${\cal M}_{GH}$ coincide components and arc components. 
Moreover, each component is open and
$\overline{\cal M}_{GH} = {\cal M}_{GH}$ is the 
topological sum of its components, 
\[ {\cal M}_{GH} = \sum\limits_{i \in I} \comp{GH} (X_i) . \]
\end{corollary}
\hfill $\Box$

\begin{proposition}
Let $X \in {\cal M}_{GH}$. Then $\comp{GH}(X)$ is given by
\[ \comp{GH}(X) = \{ Y \in {\cal M}_{GH} | d_{GH}(X,Y) < \infty \} . \]
\end{proposition}
\hfill $\Box$

We call a map $\Phi : X \lra Y$ metrically semilinear if it 
satisfies the following two conditions.

1. It is uniformly metrically proper, i. e. for each
$R>0$ there is an $S>0$ such that the inverse image under
$\Phi$ of a set of diameter $\le R$ is a set of diameter
$\le S$.

2. There exists a constant $C_\Phi \ge 0$ such that for all
$x_1, x_2 \in X$ 
$d(\Phi(x_1), \Phi(x_2)) \le d(x_1,x_2) + C_\Phi$.

Two metric spaces $X$ and $Y$ are called me\-tri\-cally semilinear
equivalent if there exist met\-ri\-cal semilinear maps
$\Phi : X \lra Y$, $\Psi : X \lra Y$ and constants $D_X$, $D_Y$ 
such that for all $x \in X$, $y \in Y$
\be
 d(x, \Psi \Phi x) \le D_X, \quad d(\Psi \Phi y, y) \le D_Y .
\ee

\begin{proposition}
$Y \in \comp{GH}(X)$, i, e. $d_{GH}(X,Y) < \infty$ if and only if
$X$ and $Y$ are metrically semilinear equivalent.
\end{proposition}

We refer to [8] for the proof.
\hfill $\Box$

At this general level there are still some important further
classes of maps in the category of proper metric spaces.

We call a map $\Phi : X \lra Y$ coarse if it is 

1. metrically proper, i. e. for each bounded set $B \subseteq Y$
the inverse image $\Phi^{-1}(B)$ is bounded in $X$, and

2. uniformly expansive, i. e. for $R>0$ there is $S>0$ s. t. 
$d(x_1, x_2) \le R$ implies $d(\Phi x_1, \Phi x_2) \le S$.

A coarse map is called rough if it is additionally uniformly
metrically proper. $X$ and $Y$ are called coarsely or roughly 
equivalent if there exist coarse or rough maps
$\Phi : X \lra Y$, $\Psi : Y \lra X$, respectively, satisfying
(2.2). A metrically proper map $\Phi : X \lra Y$ is called
Lipschitz if there holds 
$d(\Phi(x_1),\Phi(x_2)) \le C_\Phi \cdot d(x_1,x_2)$. Lipschitz
maps are continuous. $X$ and $Y$ are called coarsely Lipschitz
equivalent if there are Lipschitz maps 
$\Phi : X \lra Y$, $\Psi : Y \lra X$ satisfying (2.2).
If additionally $\Psi \Phi$ and $\Phi \Psi$ are homotopic to
$id_X$, $id_Y$ by means of a Lipschitz homotopy, respectively,
then $X$ and $Y$ are called coarsely Lipschitz homotopy
equivalent.

The following is immediately clear from the definitions.

\begin{proposition}
A metrically semilinear map is rough and a rough map is coarse.
Hence there are inclusions

metrically semilinear equivalence class of $X$ 
$(\equiv \comp{GH}(X)) \subseteq$ rough equivalence class
of $X \subseteq$ coarse equivalence class of $X$

and

coarse Lipschitz homotopy equivalence class of $X \subseteq$
coarse Lipschitz equivalence class of $X \subseteq$
coarse equivalence class of $X$.
\end{proposition}
\hfill $\Box$

For later applications to open manifolds $(M^n,g)$ we will
sharpen the Lipschitz notions by requiring not only (2.2) but
additionally controling the Lipschitz constants.

From now on we define a Lipschitz map as defined above
additionally to be uniformly metrically proper.

Define for a Lipschitz map $\Phi : X \lra Y$
\[ \dil \Phi := \sup\limits_{\begin{array}{c} x_1, x_2 \in X \\
   x_1 \neq x_2 \end{array}} 
   \frac{d(\Phi x_1, \Phi x_2)}{d(x_1,x_2)} . \]
Set
\bea
   d_L(X,Y) &:=& \inf \{ \max \{ 0, \log \dil \Phi \} +
   \max \{ 0, \log \dil \Psi \} + \sup\limits_{x \in X} 
   d(\Psi \Phi x, x) + \sup\limits_{y \in Y} d(\Phi \Psi y, y)
   \non \\
   &|& \Phi : X \lra Y, \Psi : Y \lra X \mbox{ Lipschitz maps} \} ,
   \non
\eea
if $\{ \dots \} \neq \emptyset$ and $\inf \{ \dots \}$ 
is $< \infty$ and set $d_L(X,Y) = \infty$ in the other case.
Then $d_L \ge 0$, symmetric and $d_L(X,Y) = 0$ if $X$ and
$Y$ are isometric. Set ${\cal M}_L = {\cal M} / \sim$, 
where $X \sim Y$ if $d_L(X,Y)=0$.

Let $\delta > 0$ and define
\[ V_\delta = \{ (X, Y) \in {\cal M}^2_L | d_L(X,Y) < \infty \} . \]
The proofs of the following assertions are already more
technical and lengthy. Hence we must refer to [8].

\begin{proposition}
${\cal L} = \{ V_\delta \}_{\delta>0}$ is a basis for a 
metrizable uniform structure ${\cal U}_L({\cal M}_L)$
\end{proposition}
\hfill $\Box$

Denote by ${\cal M}_L(nc)$ the (class of) noncompact proper
metric spaces.

\begin{proposition}
${\cal M}_L(nc)$ is complete with respect to 
${\cal U}_L({\cal M}_L)$.
\end{proposition}
\hfill $\Box$

\begin{proposition}
${\cal M}_L$ and $\overline{\cal M}_L(nc) = {\cal M}_L(nc)
\subset {\cal M}_L$ are locally arcwise connected.
Hence components coincide with arc components, components
are open and ${\cal M}_L$ and 
$\overline{\cal M}_L(nc) = {\cal M}_L(nc)$ have topological
sum representations
\bea
 {\cal M}_L &=& \comp{L}(point) + \sum\limits_{i \in I} \comp{L}{X_i}
 \non \\
 {\cal M}_L(nc) &=& \sum\limits_{i \in I} \comp{L}{X_i}
 \quad \mbox{and} \non \\
 \comp{L}(X) &=& \{ Y \in {\cal M}_L | d_L (X,Y) < \infty \} .
 \non
\eea
In particular all compact spaces ly in the component of the
1--point--space.
\hfill $\Box$
\end{proposition}

A sharpening of this uniform structure is given if we restrict
the maps to Lipschitz homeomphisms.

Define
\bea 
   d_{L, top} (X,Y) &=& \inf \{ \max \{ 0, \log \dil \Phi \} +
   \max \{ 0, \log \dil \Phi^{-1} \} \non \\ 
   &|& \Phi : X \lra Y \mbox{ bi--Lipschitz homeomorphism} \}
   \non
\eea
if there exists such a $\Phi$ and define $d_{L,top}(X,Y)=\infty$
if they are not bi--Lipschitz homeomorph.
Then $d_{L,top} \ge 0$, symmetric and 
$d_{L,top} (X,Y) = 0$ if $X$ and $Y$ are isometric. 
Set ${\cal M}_{L,top} = {\cal M} / \sim$, $X \sim Y$
if $d_{L,top} (X,Y) = 0$, and set
\[ V_\delta = \left\{ (X,Y) \in {\cal M}^2_{L,top} \en | \en 
   d_{L,top} (X,Y) < \delta \right\} . \]

\begin{proposition}
${\cal L} = \{ V_\delta \}_{\delta > 0}$ is a basis for a
metrizable uniform structure 
${\cal U}_{L,top}({\cal M}_{L,top})$.
\end{proposition}
\hfill $\Box$

\begin{proposition}

a) $\overline{\cal M}^{{\cal U}_{L,top}}_{L,top} = 
   {\cal M}_{L,top}$.

{\bf b)} ${\cal M}_{L,top}$ is locally arcwise connected. Hence
components coincide with arc components and components are
open.

{\bf c)} ${\cal M}_{L,top}$ has a decomposition as a topological
sum, 
\[ {\cal M}_{L,top} = \sum\limits_{i \in I} \comp{}(X_i) . \]

{\bf d)} $\comp{} (X) = \comp{L,top} (X) = \left\{ Y \in {\cal M}_{L,top}
   (X, Y)  < \infty \right\}$.
\end{proposition}
\hfill $\Box$

\no
{\bf Remark.}
We see that 2.14 a) is in ${\cal U}_{L,top}$ valid without
restriction to noncompact spaces as we due in the 
${\cal U}_L$--case.
\hfill $\Box$

Finally we define still three further uniform structures which
measure or express the homotopy neighborhoods and, secondly,
admit only compact deviations of the spaces inside one 
component.

Define
\bea 
   d_{L,h}(X,Y) &:=& \inf \Big\{ \max \{ 0, \log \dil \Phi \}
   + \max \{ 0, \log \dil \Psi \} + \sup\limits_X 
   d( \Psi \Phi x, x) + \sup\limits_Y d( \Phi \Psi y, y) 
   \non \\
   &|& \Phi : X \lra Y, \Psi : Y \lra X \mbox{ are 
   (uniformly proper) Lipschitz homotopy} \non \\
   && \en \mbox{equivalences, inverse to each other} 
   \Big\} \non
\eea
if there exist such a homotopy equivalences and set
$d_{L,h} (X,Y) = \infty$ in the other case. Here and in the
sequel we require from the homotopies to $id_X$ or $id_Y$,
respectively, that they are uniformly proper and Lipschitz.

$d_{L,h} \ge 0$, $d_{L,h}$ is symmetric and $d_{L,h}(X,Y)=0$
if $X$ and $Y$ are isometric. Define 
${\cal M}_{L,h}={\cal M} / \sim$, $X \sim Y$ if
$d_{L,h}(X,Y)=0$ and set
\[ V_\delta = \{ (X,Y) \in {\cal M}^2_{L,h} | 
   d_{L,h} (X,Y) < \delta \} . \]

\begin{proposition}
${\cal L} = \{ V_\delta \}_{\delta > 0}$ is a basis for a
metrizable uniform structure 
${\cal U}_{L,h} {\cal M}_{L,h}$.
\end{proposition}
\hfill $\Box$

\begin{proposition}

{\bf a)} $\overline{\cal M}^{{\cal U}_{L,h}}_{L,h}(nc) = 
   {\cal M}_{L,h}(nc)$.

{\bf b)} ${\cal M}_{L,h}$ is locally arcwise connected. Hence
components conincide with arc components and components
are open.

{\bf c)} ${\cal M}_{L,h}$ has a representation as a topological sum, 
\[ {\cal M}_{L,h} = \sum\limits_{i \in I} \comp{} (X_i) . \]

{\bf d)} $\comp{}(X) \equiv \comp{L,h}(X) = \{ Y \in {\cal M}_{L,h}
   | d_{L,h}(X,Y) < \infty \}$.
\hfill $\Box$
\end{proposition}

Let $\Phi: X \lra Y$, $\Psi: Y \lra X$ be Lipschitz maps. We 
say $\Phi$ and $\Psi$ are stable Lipschitz homotopy equivalences
at $\infty$ inverse to each other it there exists a compact
set $K^0_X \subset X$ s. t. for any $K^0_X \subset K_X$ there
exists $K_Y \subset Y$ s. t. $\Phi_{X \setminus K_X} : 
X \setminus K_X \lra Y \setminus K_Y$ is a Lipschitz  h. e.
with homotopy inverse $\Psi|_{Y \setminus K_Y}$ and $\Psi$
has the analogous property. 

Set
\bea 
   d_{L,h,rel}(X,Y,rel) &:=& \inf \Big\{ \max \{ 0, \log \dil \Phi \}
   + \max \{ 0, \log \dil \Psi \} + \sup\limits_X 
   d( \Psi \Phi x, x) + \sup\limits_Y d( \Phi \Psi y, y) 
   \non \\
   &|& \Phi : X \lra Y, \Psi : Y \lra X \mbox{ are 
   stable Lipschitz homotopy equivalences} \non \\
   && \en \mbox{at } \infty, \mbox{inverse to each other} 
   \Big\} \non
\eea
if there exist such a $\Phi$,  $\Psi$ and set 
$d_{L,h,rel}(X,Y) = \infty$ is the other case. Then
$d_{L,h,rel} \ge 0$, symmetric and $d_{L,h,rel} = 0$ if
$X$ and $Y$ are isometric. Set
${\cal M}_{L,h,rel}={\cal M} / \sim$, $X \sim Y$ if
$d_{L,h,rel}(X,Y)=0$ and set
\[ V_\delta = \{ (X,Y) \in {\cal M}^2_{L,h,rel} | 
   d_{L,h,rel} (X,Y) < \delta \} . \]

\begin{proposition}
${\cal L} = \{ V_\delta \}_{\delta > 0}$ is a basis for a
metrizable uniform structure 
${\cal U}_{L,h,rel}$.
\end{proposition}
\hfill $\Box$

\begin{proposition}

a) $\overline{\cal M}^{{\cal U}_{L,h,rel}}_{L,h,rel}(nc) = 
   {\cal M}_{L,h,rel}(nc)$.

{\bf b)} ${\cal M}_{L,h,rel}$ is locally arcwise connected. In 
particular components coincide with arc components, 
components are open and ${\cal M}_{L,h,rel}$ has a
representation as topological sum,  
\[ {\cal M}_{L,h,rel} = \sum\limits_{i \in I} \comp{} (X_i) , \]
where
\[ \comp{}(X) \equiv \comp{L,h,rel}(X) = \{ Y \in {\cal M}_{L,h,rel}
   | d_{L,h,rel}(X,Y) < \infty \}. \]
\end{proposition}
\hfill $\Box$

The last uniform structure ${\cal U}_{L,top,rel}$ is defined
by $d_{L,top,rel}(X,Y)$, where we require that 
$\Phi : X \lra Y$, $\Psi: Y \lra X$ are outside compact sets
bi--Lipschitz homeomorphisms, inverse to each other. We obtain
${\cal M}_{L,top,rel}$. There holds 
$\overline{\cal M}_{L,top,rel} = {\cal M}_{M,top,rel}$. The
other assumptions of 2.16 hold correspondingly.

We finish the section with a scheme which makes clear the
achievements.

\begin{center}
One coarse equivalence class

\bi
$\swarrow$ splits into $\searrow$

\bi
many GH--components \qquad \qquad many L--components,

\newpage
one L--component

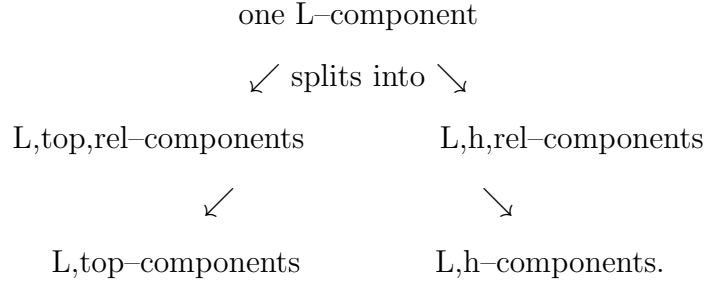
\bi
$\swarrow$ splits into $\searrow$

\bi
L,top,rel--components \qquad \qquad L,h,rel--components

\bi
$\swarrow$ \hspace{3cm} $\searrow$

\bi
L,top--components \qquad \qquad L,h--components.
\end{center}

It is now a natural observation that the classification of 
noncompact proper metric spaces splits into two main tasks

1. ''counting'' the components at any horizontal level,

2. ''counting'' the elements inside each component.

A really complete solution ot this two problems, i. e. a
complete characterization by computable and handable invariants,
is now a day hopeless. It is a similar platonic goal as the
''classification of all topological spaces''. Nevertheless
stands the task to define series of invariants which at least
permit to decide (in good cases) nonequivalence. This will be
the topic of section 5.

Finally we remark that GH--components $(d_{GH}(X,Y) < \infty)$
and L--components $(d_L(x,Y))$ are very different. Roughly 
spoken, $d_{GH}$ is in the small unsharp and in the large
relatively sharp, $d_L$ quite inverse. We refer for many
geometric examples to [8].

\section{Some materials from nonlinear global analysis on open
manifolds}

\setcounter{equation}{0}

Let $(M^n,g)$ be a Riemannian manifold. We consider the conditions (I) and
(B$_k$),\\
\begin{tabular}{cl}
(I) &  $r_{inj}(M^n,g)=\inf\limits_{x\in M} r_{inj}(x) >0$, \\
(B$_k$) \qquad & $|\nabla^iR|\le C_i$, $0\le i \le k$,  
\end{tabular}\\
where $r_{inj}$ denotes the injectivity radius, $\nabla=\nabla^g$ the
Levi-Civita connection, $R=R^g$ the curvature and $|\cdot|$ the pointwise
norm. $(M^n,g)$ has bounded geometry of order $k$ if it satisfies the conditions
(I) and (B$_k$). Every compact manifold $(M^n,g)$ or homogeneous Riemannian
space or Riemannian covering $(\tilde{M},\tilde{g})$ of a compact manifold
$(M^n,g)$ satisfies (I) and (B$_\infty$). More general, given $M^n$ open and
$0\le k\le \infty$, there exists a complete metric of bounded geometry of
order $k$ (cf. [19]), i.e. the existence of such a metric does not
restrict the underlying topological type. (I) implies completeness. Let
$(E,h) \lra (M^n,g)$ be a Riemannian vector bundle and $\nabla=\nabla^h$ a metric
connection with respect to $h$. Quite analogously we consider the condition\\
\begin{tabular}{cl}
(B$_k(E,\nabla)$) \qquad & $|\nabla^iR^E|\le C_i$, \qquad $0\le i \le k$,  
\end{tabular}\\
where $R^E$ denotes the curvature of $(E,\nabla)$.
\begin{lemma}\label{l21}
  If $(M^n,g)$ satisfies $(B_k)$ and ${\cal U}$ is an atlas of normal coordinate
  charts of radius $\le r_0$ then there exist constants $C_\alpha$,
  $C'_\beta$ such that
  \begin{eqnarray}\label{e21}
    && |D^\alpha g_{ij}| \le C_\alpha, \quad |\alpha|\le k, \\
    \label{e22}
    && |D^\beta \Gamma_{ij}^m| \le C'_\beta, \quad |\beta|\le k-1,
  \end{eqnarray}
  where $C_\alpha$, $C'_\beta$ are independent of the base points of the
  normal charts and depend only on $r_0$ and on curvature bounds including
  bounds for the derivations.
\end{lemma}
We refer to [11] for the rather long and technical proof which uses
iterated inhomogeneous Jacobi equations.
\hfill $\Box$

This lemma carries over to the case of Riemannian vector bundles. A normal
chart $U$ in $M$ of radius $\le r_0$, an orthonormal frame $e_1,\dots,e_N$
over the base point $p\in U\subset M$ and its radial parallel translation
define a local orthonormal frame field $e_\alpha$, a so called synchronous
frame and connection coefficients $\Gamma_{\alpha i}^\beta$ by
\[ \nabla_{\frac{\partial}{\partial x_i}} e_\alpha = \Gamma_{\alpha i}^\beta
e_\beta. \] 
\begin{lemma}\label{l22}
  Assume $(B_k(M^n,g))$, $(B_k(E,\nabla))$, $k\ge 1$, and $\Gamma_{\alpha
    i}^\beta$ as above. Then
  \begin{eqnarray}
    |D^\gamma\Gamma_{\alpha i}^\beta| \le C_\gamma, \quad |\gamma|\le k-1,
    \quad \alpha,\beta=1,\dots,N,\quad i=1,\dots,n,
  \end{eqnarray}
  where the $C_\gamma$ are constants depending on curvature bounds, $r_0$, and
  are independent of $U$.
\end{lemma}
We refer to [11] for the proof.
\hfill $\Box$

We recall for what follows some simple facts concerning Sobolev spaces on open
manifolds. Let $(E,h,\nabla^h) \lra (M^n,g)$ be a Riemannian vector bundle. Then
the Levi-Civita connection $\nabla^g$ and $\nabla^h$ define metric connections
$\nabla$ in all tensor bundles $T_v^u\otimes E$. Denote by
$C^\infty(T_v^u\otimes E)$ all smooth sections, $C^\infty_c(T_v^u\otimes E)$
those with compact support. In the sequel we write $E$ instead of
$T_v^u\otimes E$, keeping in mind that $E$ can be an arbitrary vector
bundle. Now we define for $p\in \R$, $1\le p <\infty$ and $r$ a nonnegative
integer
\bea
  |\varphi|_{p,r} &:= & \left(\int\sum\limits_{i=0}^r |\nabla^i\varphi|^p_x
    dvol_x(g)\right)^{1/p} \non \\
  \Omega^p_r(E) & = & \{ \varphi\in C^\infty(E) \mid |\varphi|_{p,r} < \infty
  \}, \non \\
  \overline{\Omega}^{p,r}(E) & = & \mbox{completion of } \Omega^p_r(E) 
   \mbox{ with respect  to} |\cdot |_{p,r}, \non \\
  \begin{array}{c} {\rm _o} \\[-0,5ex] \Omega \\[-0,5ex] {} \end{array}
  {}^{p,r}(E) 
  & = & \mbox{completion of } C^\infty_c(E) \mbox{ with respect
    to } |\cdot |_{p,r} \mbox{ and} \non \\
  \Omega^{p,r}(E) & = & \{\varphi\mid\varphi \mbox{ measurable distributional
    section with } |\varphi|_{p,r}<\infty \}. \non
\eea
Furthermore, we define
\bea
  ^{b,m}|\varphi| & := & \sum\limits_{i=0}^m\sup\limits_x|\nabla^i\varphi|_x, \non \\
  ^{b,m}\Omega(E) & = & \{ \varphi\mid\varphi\enspace C^m- \mbox{section and}
   \enspace ^{b,m}|\varphi|<\infty \} \mbox{ and} \non \\
  ^{b,m}\begin{array}{c} {\rm _o} \\[-0,5ex] \Omega \\[-0,5ex] {} \end{array}
   (E) & = & \mbox{completion of } C^\infty_c(E) 
   \mbox{ with respect to }
    ^{b,m}|\cdot|. \non 
\eea
$^{b,m}\Omega(E)$ equals the completion of
\[ ^b_m\Omega(E) = \{\varphi\in C^\infty(E) \mid ~^{b,m}|\varphi|<\infty\} \]
with respect to $^{b,m}|\cdot|$.
\begin{proposition}
  The spaces $\begin{array}{c} {\rm _o} \\[-0,5ex] \Omega \\[-0,5ex] {} 
  \end{array}{}^{p,r}
  (E)$, $\bar{\Omega}^{p,r}(E)$, $\Omega^{p,r}(E)$,
  $^{b,m}\begin{array}{c} {\rm _o} \\[-0,5ex] \Omega \\[-0,5ex] {} \end{array}
  (E)$, $^{b,m}\Omega(E)$ are Banach spaces and there are
  inclusions
  \begin{eqnarray*}
    & \begin{array}{c} {\rm _o} \\[-0,5ex] \Omega \\[-0,5ex] {} \end{array}
    {}^{p,r}(E)\subseteq \bar{\Omega}^{p,r}(E) \subseteq
    \Omega^{p,r}(E),\\
    & ^{b,m}\begin{array}{c} {\rm _o} \\[-0,5ex] \Omega \\[-0,5ex] {} \end{array}
   (E)\subseteq ~^{b,m}\Omega(E).
  \end{eqnarray*}
  If $p=2$ then $\begin{array}{c} {\rm _o} \\[-0,5ex] \Omega \\[-0,5ex] {} \end{array}{}
  ^{2,r}(E)$, $\bar{\Omega}^{2,r}(E)$, $\Omega^{2,r}(E)$
  are Hilbert spaces.
\hfill $\Box$
\end{proposition}
$\begin{array}{c} {\rm _o} \\[-0,5ex] \Omega \\[-0,5ex] {} \end{array}
{}^{p,r}(E)$, $\bar{\Omega}^{p,r}(E)$, $\Omega^{p,r}(E)$ are different in
general.
\begin{proposition}
  If $(M^n,g)$ satisfies $(I)$ and $(B_k)$ then
  \[ \begin{array}{c} {\rm _o} \\[-0,5ex] \Omega \\[-0,5ex] {} \end{array}
  {}^{p,r}(E) = \bar{\Omega}^{p,r}(E) = \Omega^{p,r}(E), \quad 0 \le r
  \le k+2. \]
\end{proposition}
We refer to [12] for the proof.
\hfill $\Box$

\begin{theorem}
  Let $(E,h,\nabla^E) \lra (M^n,g)$ be a Riemannian vector bundle satisfying $(I)$,
  $(B_k$$(M^n,g))$, $(B_k(E,\nabla))$, $k\ge 1$.
  \begin{description}
  \item[{\bf a)}] Assume $k\ge r$, $r-\frac{n}{p}\ge s-\frac{n}{q}$, $r\ge s$, $q\ge
    p$. Then
    \begin{equation}\label{e33}
      \Omega^{p,r}(E) \hookrightarrow \Omega^{q,s}(E)
    \end{equation}
    continuously.
  \item[{\bf b)}] If $r-\frac{n}{p} > s$, then
    \begin{equation}
      \Omega^{p,r}(E) \hookrightarrow ~^{b,s}\Omega(E)
    \end{equation}
    continuously.
  \end{description}
\end{theorem}
\hfill $\Box$

A key role for many calculations and estimates in nonlinear
global analysis plays the module structure theorem which asserts under
which conditions the tensor product of two Sobolev sections is again
a Sobolev section. We do not use this theorem here explicitely and 
refer to [13].

We consider in section 4 spaces of metrics and give now a very small
review of that. Let $M^n$ be an open smooth manifold, 
${\cal M} = {\cal M}(M)$ the space of all complete Riemannian metrics.
Let $g \in {\cal M}$. We define
\[ {}^bU(g) = \{ g' \in {\cal M} | ^b|g-g'|_{g} := \sup\limits_{x \in M}
   |g-g'|_{g,x} < \infty, ^b|g-g'|_{g'} < \infty  \} . \]
It is easy to see that $^bU(g)$ coincides with the quasi isometry
class of $g$, i. e.  $g' \in ^bU(g)$ if and only if there exist
$c,c' >0$ such that
\[ c \cdot g' \le g \le c' \cdot g' . \]
Denote for $g, g' \in {\cal M}$ by $\nabla = \nabla^g$,
$\nabla' = \nabla^{g'}$ the Levi-Civita connections. Set for
$m \ge 1$, $\delta > 0$
\bea 
    V_\delta &=& \{ (g,g') \in {\cal M}^2 | ^b|g-g'|_g <\delta, 
   ^b|g-g'|_{g'} <\delta \mbox{ and } \non \\
   ^{b,m}|g-g'|_g &:=&
   ^b|g-g'|_g + \sum\limits^{m-1}_{j=0}
   {}^b|\nabla^j(\nabla-\nabla')|_g <\delta   \} . \non
\eea

\begin{proposition}
The set ${\cal B} = \{ V_\delta \}_{\delta>0}$ is a basis ${\cal B}$
for a metrizable uniform structure on ${\cal M}$.
\end{proposition}

Denote by $^{b,m}{\cal M}$ the corresponding completed uniform
space.

\begin{proposition}
The space $^{b,m}{\cal M}$ is locally contractible.
\end{proposition}

\begin{corollary}
In $^{b,m}{\cal M}$ components and arc components coincide.
\end{corollary}

Set
\[ ^{b,m}U(g) = \{ g' \in ^{b,m}{\cal M} | ^b|g-g'|_{g'} < \infty, 
    ^{b,m}|g-g'|_g < \infty \} . \]

\begin{proposition}
Denote by $\comp{}(g)$ the component of $g \in ^{b,m}{\cal M}$. Then
\[ \comp{}(g) = ^{b,m} U(g). \]
\end{proposition}

\begin{theorem}
The space $^{b,m}{\cal M}$ has a representation as a topological sum
\[ ^{b,m}{\cal M} = \sum\limits_{i \in I} ^{b,m} U(g_i). \]
\end{theorem}

\begin{theorem}
Each component of $^{b,m}{\cal M}$ is a Banach manifold.
\end{theorem}

Denote for given $M$
\[ {\cal M} (I,B_k)  = \{ g \in {\cal M} | g \mbox{ satisfies } 
   (I) \mbox{ and } (B_k) \} . \]

\me
\no
{\bf Remark.} 
$(I)$ always implies completeness.
\hfill $\Box$

Metrics of bounded geometry wear a natural inner Sobolev topology.
Let $k \ge r < \frac{n}{p} + 1$, $\delta > 0$ and set
\bea 
  V_\delta &=& \{ (g,g') \in {\cal M} (I, B_k)^2 | ^b|g-g'|_g < \delta,
  ^b|g-g'|_{g'} < \delta \mbox{ and} \non \\
  |g-g'|_{q,p,r} &:=& ( \int ( |g-g'|_{g,x} + \sum\limits^{r-1}_{i=0}
  |\nabla^i (\nabla-\nabla' )|^p_{g,x} dvol_x(g) )^{\frac{1}{p}} 
  < \delta \} . \non
\eea
 
\begin{proposition}
The set $\{ V_\delta \}_{\delta>0}$ is a basis 
for a metrizable uniform structure on ${\cal M}(I,B_k)$.
\end{proposition}

Denote by ${\cal M}^{p,r}$ the corresponding completed uniform
space.

\begin{proposition}
The space ${\cal M}^{p,r}(I,B_k)$ is locally contractible.
\end{proposition}

\begin{corollary}
In ${\cal M}^{p,r}(I,B_k)$ components and arc components coincide.
\end{corollary}

Set for $g \in {\cal M}^{p,r}(I,B_k)$
\[ U^{p,r}(g) = \{ g' \in {\cal M}^{p,r}(I,B_k) | ^b|g-g'|_g < \infty, 
   ^b|g-g'|_{g'} < \infty, ^b|g-g'|_{g,p,r} < \infty \} .  \]

\begin{proposition}
Denote by $\comp{}(g)$ the component of $g \in {\cal M}^{p,r}(I,B_k)$. Then
\[ \comp{}(g) = U^{p,r} (g). \]
\end{proposition}

\begin{theorem}
Let $M^n$ be open, $k \ge r > \frac{n}{p}+1$. Then ${\cal M}^{p,r}(I,B_k)$
has a representation as a topological sum
\[ {\cal M}^{p,r}(I,B_k) = \sum\limits_{i \in I} U^{p,r}(g_i) . \] 
\end{theorem}

We refer to [14] for all proofs.

Finally we give one hint to the manifolds of maps theory. Let $(M^n,g)$,
$(N^{n'},h)$ be open, satisfying $(I)$ and $(B_k)$ and let
$f \in C^\infty(M,N)$. Then the differential $df = f_* = Tf$ is a section
of $T^*M \otimes f^*TN$. $f^*TN$ is endowed with the induced connection
$f^* \nabla^h$. The connections $\nabla^g$ and $f^*\nabla^h$ induce
connections $\nabla$ in all tensor bundles
$T^q_s(M) \otimes f^*T^u_v(N)$. Therefore, $\nabla^mdf$ is well defined.
Assume $m \le k$. We denote by $C^{\infty,m}(M,N)$ the set of all
$f \in C^{\infty}(M,N)$ satisfying
\[ ^{b,m} |df| = \sum\limits^{m-1}_{i=0} \sup\limits_{x \in M}
    |\nabla^i df|_x < \infty . \]
It is now possible to define for $C^{\infty,m}(M,N)$ uniform structures
to obtain manifolds of maps $^{b,m}\Omega(M,N)$, $\Omega^{p,r}(M,N)$,
manifolds of diffeomorphisms $^{b,m}D(M,N)$, $D^{p,r}(M,N)$ and groups
of diffeomorphisms $^{b,m}D(M)$, $D^{p,r}(M)$.

But this approach is extraordinarily complicated and extensive. We
refer to [15].

\me
\no
{\bf Remark.}
For closed manifolds all these things are very simple.
\hfill $\Box$

\section{Uniform structures of open manifolds}

\setcounter{equation}{0}

Any complete Riemannian manifold $(M^n,g)$ defines a proper
metric space and hence an element of ${\cal M}$, ${\cal M}_{GH}$, 
${\cal M}_{L}$, ${\cal M}_{L,top,rel}$, ${\cal M}_{L,top}$, 
${\cal M}_{L,h,rel}$, ${\cal M}_{L,h}$. Denote by ${\cal M}^n(mf)$
the subset of (classes of) complete Riemannian $n$--manifolds
$(M^n,g)$. The restriction of an uniform structure to a subset
yields an uniform structure and we obtain uniform spaces
$({\cal M}^n_{GH}(mf), {\cal U}_{GH} |_{{\cal M}^n_{GH}}(mf))$,
$({\cal M}^n_L(mf), {\cal U}_L |_{{\cal M}^n_L}(mf))$ etc.. 
These uniform structures do not take into account the smooth and
Riemannian structure but only the (distance) metrical structure.
Metrically semilinear and Lipschitz maps which enter into the
definition of $d_{GH}$, $d_L$ etc. are far from being smooth.
On the other hand, we defined in [14] and in the last section
for one fixed manifold
$M^n$ several uniform structures of Riemannian metrics. Roughly
spoken, the space ${\cal M}(M)$ of complete Riemannian metrics
$g$ on $M$ splits into a topological sum
${\cal M}(M) = \sum\limits_{i \in I} \comp{}(g_i)$ depending on the
norm in question. Then the fundamental question arises how are
$\comp{}(M,g_i)$ and $\comp{GH}(M,g_i)$, $\comp{L}(M,g_i)$ etc.
related? We will give very shortly some partial answers in this
section. Moreover, we will generalize the (analytically defined)
uniform structures of metrics for one fixed manifold to uniform
structures of Riemannian manifolds. Still other questions concern
the completions of ${\cal M}^n(mf)$ with respect to the
considered uniform structures and the completions of certain
subspaces of metrics, e. g. Ricci curvature $\le 0$. But these
more differential geometric questions are outside the space of
this paper and will be studied in [16]. Looking at our general
approach, the first main interesting questions are the relations
between the components defined until now.

\begin{proposition}
Let $(M^n,g) \in {\cal M}^n(mf)$, 
$g' \in {}^{b,m} \comp{}(g) \subset {}^{b,m} {\cal M}$,
$m \ge 0$. Then there holds

{\bf a)} $(M,g') \in \comp{L,top} (M,g)$,

{\bf b)} $(M,g') \in \comp{L} (M,g)$
\end{proposition}

We refer to [8] for the proof.
\hfill $\Box$

\no
{\bf Remark.}
We cannot prove this for $d_{GH}$. $d_{GH}$ is locally very rough
but measures the metric relations in the large relatively exact.
But this property does not immediately follow from 
$g' \in {}^{b,m} \comp{}(g)$.
\hfill $\Box$

\begin{corollary}
The assertions a) and b) hold if
$g' \in \comp{}^{p,r}(g) \subset {\cal M}^{p,r}(I,B_k)$,
$k \ge r > \frac{n}{p}+1$.
\end{corollary}
\hfill $\Box$

To admit the variation of $M$ in $(M^n,g)$ too, we define still
another uniform structure. First we admit arbitrary complete
metrics, i. e. we do not restrict to metrics of bounded geometry.
Consider complete manifolds $(M^n,g)$ $(M'^n,g')$ and 
$C^{\infty,m}(M,M')$. A diffeomorphism $f : M \lra M'$ will be
called $m$--bibounded if $f \in C^{\infty,m}(M,M')$ and
$f^{-1} \in C^{\infty,m}(M',M)$. Sufficient for this is
1. $f$ is a diffeomorphism, 2. $f \in C^{\infty,m}(M,M')$,
3. $\inf\limits_x | \lambda |_{\min} (df)_x > 0$.

Let $\delta > 0$ and set
\bea 
  V_\delta &=& \{ ((M^n_1,g_1), (M^n_2,g_2)) \in {\cal M}^n(mf)^2 
  \en | \en \mbox{There exists a diffeomorphism } \non \\
  && f: M_1 \lra M_2, \, f \, (m+1)-\mbox{bibounded} \non \\
  && (1+\delta+\delta\sqrt{2n(n-1)})^{-1} \cdot g_1 \le f^*g_2 \le
  (1+\delta+\delta\sqrt{2n(n-1)}) \cdot g_1  \} . \non
\eea

\begin{proposition}
${\cal L} = \{ V_\delta \}_{\delta > 0}$  is a basis for
a uniform structure
${}^{b,m}{\cal U}_{diff}({\cal M}^m(mf))$
\end{proposition}

We refer to [8] for the long and technical proof.
\hfill $\Box$

Denote ${\cal M}^n(mf)$ endowed with the 
${}^{b,m}{\cal U}_{diff}$--topology by 
${}^{b,m}{\cal M}^n(mf)$.

\begin{proposition}
${}^{b,m}{\cal M}^n(mf)$ is locally arcwise connected.
\hfill $\Box$
\end{proposition}

\begin{corollary}
In ${}^{b,m}{\cal M}^n(mf)$ components coincide with arc
components.
\end{corollary}
\hfill $\Box$

\begin{theorem}

{\bf a)} ${}^{b,m}{\cal M}^n(mf)$ has a representation as topological
sum,
\[ {}^{b,m}{\cal M}^n(mf) = \sum\limits_i 
   {}^{b,m} \comp{diff}(M_i,g_i) . \]

{\bf b)} 
\bea 
  {}^{b,m} \comp{diff}(M,g) &=& \{ (M',g') \in {\cal M}^n(mf) 
  \en | \en
  \mbox{There exists a diffeomorphism } \non \\
  &&   f : M \lra M' \in  C^{\infty,m+1}(M,M') 
   \mbox{ s. t. } {}^{b,m} |f^*g'-g|_g < \infty \} . \non 
\eea
\end{theorem}

We refer to [8] for the proof.
\hfill $\Box$

\no
{\bf Remark.}
If $f : (M_0,g_0) \lra (M_1,g_1)$ is a diffeomorphism such that
$c_1 \cdot g_0 \le f^* g_1 \le c_2 \cdot g_0$ then
$\dil (f) \le c_2$ and $\dil (f^{-1}) \le \frac{1}{c_1}$.
\hfill $\Box$

\begin{corollary}
If $(M_1,g_1) \in {}^{b,m} \comp{diff}(M_0,g_0)$ then 
\[ d_{L,top}((M_0,g_0),(M_1,g_1)), d_L((M_0,g_0),(M_1,g_1))<\infty.\]
\end{corollary}
\hfill $\Box$

\no
{\bf Remark.}
If we assume $(I)$, $(B_k)$, $k \ge m$, then we can complete the
$(m+1)$--bibounded diffeomorphisms of $C^{\infty,m+1}(M_0,M_1)$ to get
${}^{b,m+1} \diff (M_0,M_1)$ and 4.3 -- 4.7 remain valid with 
$C^{m+1}$--diffeomorphisms (bibounded) between manifolds of bounded
geometry.
\hfill $\Box$

Now we define the uniform structures for the Riemannian case which
are parallel to them defined at the end of section 2.

Consider pairs $(M^n_1,g_1), (M^n_2,g_2) \in {\cal M}^n(mf)$
with the following property.
\bea
  && \mbox{There exist compact submanifolds } K^n_1 \subset M^n_1,
  K^n_2 \subset M^n_2 \non \\
  && \mbox{ and an isometry } 
  \Phi : M_1 \setminus K_1 \lra M_2 \setminus K_2 . 
\eea
For such pairs define
\bea 
  {}^bd_{L,iso,rel} ((M_1,g_1),(M_2,g_2)) &:=& 
  \inf \{ \max \{ 0, \log {}^b |df| \} + \max \{ 0, \log {}^b |dh| \}
  \non \\ 
  &+& \sup\limits_{x \in M_1} \dist (x,hfx) 
  + \sup\limits_{y \in M_2} \dist (y,fhy) \non \\
  &|& f \in C^{\infty} (M_1,M_2), g \in C^{\infty} (M_2,M_1)
  \mbox{ and for some } \non \\
  && K_1 \subset K \mbox{ holds } 
  f|_{M_1 \setminus K_1} \mbox{ is an isometry and } \non \\
  &&  g|_{f(M_1 \setminus K)} = f^{-1} \} . \non
\eea
If $(M_1,g_1)$, $(M_2,g_2)$ satisfy (4.1) then 
$\{ \dots \} \neq \emptyset$ and 
${}^bd_{L,iso,rel} (M_1,M_2) = \inf \{ \dots \} < \infty$. If
$(M_1,g_1)$, $(M_2,g_2)$ do not satisfy (4.1) then we define
${}^bd_{L,iso,rel} ((M_1,g_1),(M_2,g_2)) = \infty$.
${}^bd_{L,iso,rel} (\cdot , \cdot)$ is $\ge 0$, symmetric and
${}^bd_{L,iso,rel} ((M_1,g_1),(M_2,g_2)) \le \infty$.
${}^bd_{L,iso,rel} ((M_1,g_1),(M_2,g_2)) = 0$ if $(M_1,g_1)$ and
$(M_2,g_2)$ are isometric.

\me
\no
{\bf Remarks.}

\no
1) The notion Riemannian isometry and distance isometry coincide
for Riemannian manifolds. Moreover, for an isometry $f$ holds
${}^b |df|=1$.

\no
2) Any $f$ which enters into the definition of $d_{L,iso,rel}$ is
automatically an element of $C^{\infty,m}(M_1,M_2)$ for all $m$.
The same holds for $g$.
\hfill $\Box$

We denote ${\cal M}^n_{L,iso,rel} (mf) = {\cal M}^n (mf) / \sim$
where $(M_1,g_1) \sim (M_2,g_2)$ if 
\[ {}^bd_{L,iso,rel} ((M_1,g_1),(M_2,g_2)) = 0. \]

Set
\[ V_\delta = \{ ((M_1,g_1),(M_2,g_2)) \in ({\cal M}^n_{L,iso,rel}
   (mf))^2 \en | \en {}^bd_{L,iso,rel} ((M_1,g_1),(M_2,g_2)) < \delta 
   \} . \]

\begin{proposition}
${\cal L} = \{ V_\delta \}_{\delta > 0}$ is a basis for a metrizable
uniform structure ${}^b{\cal U}_{L,iso,rel}$.
\end{proposition}
\hfill $\Box$

Denote by ${}^b{\cal M}^n_{L,iso,rel}(mf)$ the corresponding uniform
space.

\begin{proposition}
If $r_{inj}(M_i,g_i) = r_i > 0$, $r = \min \{ r_1, r_2 \}$ and
${}^bd_{L,iso,rel} ((M_1,g_1),(M_2,g_2)) < r$ then $M_1$, $M_2$
are (uniformly proper) bi--Lipschitz homotopy equivalent.
\end{proposition}
\hfill $\Box$

\begin{corollary}
If we restrict to open manifolds with injectivity radius $\ge r$
then manifolds $(M_1,g_1)$, $(M_2,g_2)$ with
${}^bd_{L,iso,rel}$--distance $< r$ are automatically (uniformly
proper) bi--Lipschitz homotopy equivalent.
\end{corollary}
\hfill $\Box$

\no
{\bf Remark.}
If $(M_1,g_1)$ satisfies $(I)$ or $(I)$ and $(B_k)$ and
${}^bd_{L,iso,rel} (M_1,g_1), (M_2,g_2)) < \infty$ then
$(M_2,g_2)$ also satisfies $(I)$ or $(I)$ and $(B_k)$.
\hfill $\Box$

We cannot show that ${}^b{\cal M}^n_{L,iso,rel}$ is locally
arcwise connected, that components coincide with arc components and
that ${}^b \comp{L,iso,rel} (M,g) = \{ (M',g') | {}^bd_{L,iso,rel}
((M,g),(M',g')) < \infty \}$. The background for this is the fact
that it is impossible to connect non homotopy equivalent manifolds
by a continuous family of manifolds. A parametrization of nontrivial
surgery always contains bifurcation levels where we leave the
category of manifolds. A very simple handable case comes from
4.10.

\begin{corollary}
If we restrict ${}^b{\cal U}_{L,iso,rel}$ to open manifolds with
injectivity radius $\ge r > 0$ then the manifolds in each arc
component of this subspace are bi--Lipschitz homotopy equivalent.
\end{corollary}

\me
\no
{\bf Proof.}
This subspace is locally arcwise connected, components coincide
with arc components. Consider an (arc) component, elements $(M_1,g_1)$,
$(M_2,g_2)$ of it, connect them by an arc, cover this arc by
sufficiently small balls and apply 4.10.
\hfill $\Box$

It follows immediately from the definition that
${}^bd_{L,iso,rel}((M_1,g_1), (M_2,g_2)) <\infty$ implies
$d_L((M_1,g_1), (M_2,g_2)) <\infty$. Hence
$(M_2,g_2) \in \comp{L}(M_1,g_1)$, i. e.
\be
 \{ (M_2,g_2) \in {\cal M}^n(mf) | {}^bd_{L,iso,rel} 
 ((M_1,g_1),(M_2,g_2)) < \infty \} \subseteq \comp{L} (M_1,g_1).
\ee
For this reason we denote the left hand side $\{ \dots \}$ of
(4.2) by ${}^b \comp{L,iso,rel}(M_1,g_1) = \{ \dots \} = 
\{ \dots \} \cap \comp{L}(M_1,g_1)$ keeping in mind that this is
not an arc component but a subset (of manifolds) of a Lipschitz
arc component.

If one fixes $(M_1,g_1)$ then one has in special cases a good
overview on the elements of ${}^b \comp{L,iso,rel} (M_1,g_1)$.

\me
\no
{\bf Example.}
Let $(M_1,g_1)=(\R^n, g_{standard})$. Then
${}^b \comp{L,iso,rel} (\R^n,g_{standard})$ is in 
1--1--relation to $\{ (M^n,g) | M^n$ is a closed manifold, $g$
is arbitrary but flat in an annulus contained in disc 
neighborhood of a point$\}$.
\hfill $\Box$

This can be generalized as follows.

\begin{theorem}
Any component ${}^b \comp{L,iso,rel} (M,g)$ contains at most
countably many diffeomorphism types.
\end{theorem}

\me
\no
{\bf Proof.}
Fix $(M,g) \in {}^b \comp{L,iso,rel} (M,g)$, an exhaustion
$K_1 \subset K_2 \subset \dots$, $\bigcup K_i =M$, by
compact submanifolds and consider 
$(M',g') \in {}^b \comp{L,iso,rel} (M,g)$. Then there exist
$K' \subset M'$ and $K_i \subset M$ such that $M \setminus K_i$
and $M' \setminus K'$ are isometric. The diffeomorphism type
of $M'$ is completely determined by diffeomorphism type of
the pair $(K_1 \bigcup\limits_{\partial K_1 \cong \partial K'}
K', K_1)$ but there are only at most countably many types of
such pairs (after fixing $M$ and $K_1 \subset K_2 \subset \dots$).
\hfill $\Box$

That is, after fixing $(M,g)$, the diffeomorphism classification
of the elements in ${}^b \comp{L,iso,rel}$ $(M,g)$ seems to be
reduced to a ''handable'' countable discrete problem. This is
in fact the case in a sense which is parallel to the
classification of compact manifolds. This will be carefully
and detailed discussed in [8] -- [10]. A key role for this
plays the following

\begin{lemma}
Let $\sigma_e (\Delta_q)$ be the essential spectrum of the
Laplace operator acting on $q$--forms and assume 
$(M',g') \in {}^b \comp{L,iso,rel} (M,g)$. Then
$\sigma_e (\Delta_q) (M,g) = \sigma_e (\Delta_q) (M',g')$,
$0 \le q \le n$.
\end{lemma}
\hfill $\Box$

A further step in this classification approach will be the
definition of characteristic numbers for pairs $(M,M')$ and
of bordism which will be the content of section 5.

Instead of requiring isometry at infinity we can focus our
attention to homotopy porperties and define
\bea
  {}^{b,m}d_{L,h} ((M_1,g_1),(M_2,g_2)) &=&
  \inf \{ \max \{ 0, \log {}^b |df| \} + \max \{ 0,
  \log {}^b |dh| \}  \non \\
  &+&  \sup\limits_{x \in M_1} \dist (x,hfx)
  + \sup\limits_{y \in M_2} \dist (y,fhy) \non \\ 
  &|& f \in C^{\infty,m} (M_1,M_2), h \in C^{\infty,m} (M_2,M_1),
  f \mbox{ and } h \mbox{ are inverse} \non \\
  && \mbox{to each other uniformly
  proper homotopy equivalences} \} ,
\eea
if $\{ \dots \} \neq \emptyset$ and if $\inf \{ \dots \} < \infty$.
In the other case define ${}^{b,m}d_{L,h}((M_1,g_1),(M_2,g_2))
= \infty$. Then $d_{L,h} \ge 0$, symmetric and $=0$ if
$(M_1,g_1)$ ,$(M_2,g_2)$ are isometric. Define
${\cal M}^n_{L,h}(mf) = {\cal M}^n(mf) / \sim$ where $\sim$
means ${}^{b,m}d_{L,h}$--distance $=0$.

Set
\[ V_\delta = \{ ((M_1,g_1),(M_2,g_2)) \in ({\cal M}^n_{L,h}
   (mf))^2 \en | \en {}^{b,m}d_{L,h} ((M_1,g_1),(M_2,g_2))
   < \delta\} . \]

\begin{proposition}
${\cal L} = \{ V_\delta \}_{\delta >0}$ is a basis for a
metrizable uniform structure ${}^m{\cal U}_{L,h}$.
\end{proposition}
\hfill $\Box$

Denote by ${}^{b,m}{\cal M}_{L,h}$ the corresponding uniform
space. Here again we cannot prove that 
${}^{b,m}{\cal M}^n_{L,h}(mf)$ is locally arcwise connected.
Nevertheless, ${}^{b,m}d_{L,h} ((M_1,g_1),(M_2,g_2)) < \infty$
implies $d_L ((M_1,g_1),(M_2,g_2)) < \infty$ and hence 
$(M_2,g_2) \in \comp{L}(M_1,g_1)$.
$\{ (M_2,g_2) \in {\cal M}^n_{L,h} (mf)$ 
$| {}^{b,m}d_{L,h} ((M_1,g_1),(M_2,g_2)) < \infty \} \subseteq
\comp{L}(M_1,g_1)$, and we write ${}^{b,m} \comp{L,h}(M_1,g_1)$
for this set.

\me
\no
{\bf Remark.}
${}^{b,m} d_{L,h} ((M_1,g_1),(M_2,g_2)) = 0$ implies
$d_{L,h} ((M_1,g_1),(M_2,g_2)) = 0$. The corresponding 
implication holds in all preceding cases.
\hfill $\Box$

In a quite analogous manner as in section 2, proposition 2.15, 
and as above, i. e. restricting to maps $\in C^{\infty,m}$,
we can define ${}^{b,m} d_{L,h,rel}$, 
${}^{b,m} {\cal M}_{L,h,rel}(mf)$ and
${}^{b,m} \comp{L,h,rel}(M,g)$, where
${}^{b,m} \comp{L,h,rel}(M,g) = \{ (M',g') | 
{}^{b,m} d_{L,h,rel} ((M,g),(M',g')) < \infty \} \subset
\comp{L,h,rel}(M,g)$.

All uniform structures defined until now for manifolds are
based on the Banach ${}^{b,m} |\,\,|$--theory. It is
possible to construct an extensive theory of uniform
structures of manifolds based on the $|\,\,|_{p,r}$--Sobolev
approach. For spaces of metrics and manifolds of maps
this has been done e.g. in [14], [15]. We cannot discuss
here the Sobolev approach for reasons of space and refer
to [8]. An important philosophical hint shall be given.
The ${}^{b,m} |\,\,|$--approach for manifolds is more
related to the $d_L$--uniform structures (as we pointed out),
part of the Sobolev approach for manifolds is more related
to $d_{GH}$--uniform structures.

We conclude this section with the hint to [8] where most
of the very long details are represented.

\section{Bordism groups for open manifolds}

\setcounter{equation}{0}

We sketch very shortly our approach to bordism theory for open
manifolds. Let $(M^n,g)$, $(M'^n,g')$ be open, oriented, 
complete. We say $(M^n,g)$ is bordant to $(M'^n,g')$ if
there exists an oriented complete manifold $(B^{n+1},g_B)$ with 
boundary $\partial B$ such that the following holds.

1) $(\partial B, g_B |_{\partial B}) = (M,g) \cup (-M',g')$.
Here $=$ stands for isometry.

2) There exists a uniform Riemannian collar
$\Phi : \partial B \times [ 0, \delta [
 \begin{array}{c} \cong \\[-2ex] \lra \\[-2ex] {} \end{array}
 {\cal U}_\delta (\partial B) \subset B$, 
$\Phi^* (g_B |_{{\cal U}_\delta (\partial B)} ) =
 g_{\partial B} + dt^2$.

3) There exists $R>0$ such that $B \subseteq {\cal U}_R (M)$,
$B \subseteq {\cal U}_R (M')$.

\me
\no
{\bf Remark.}
Condition 3) looks like $d^B_H (M,M') \le R$, 
$d_{GH} (M,M') \le R$. But this is not necessary the case since 
1) and 2) do not imply that $\partial B = M \cup M'$ is
isometrically embedded as metric length space into the metric
length space $B$. $\partial B$ is isometrically embedded as
Riemannian manifold but its inner length metric will not be
the induced length metric from $B$, even not if $\partial B$
is totally geodesic as we assume by 2).
\hfill $\Box$

We denote $(M,g) 
\begin{array}{c} {} \\[-1ex] \sim \\[-1ex] b \end{array}
(M',g')$. $B^{n+1}$ is called a bordism.

\begin{lemma}
$\begin{array}{c} {} \\[-1ex] \sim \\[-1ex] b \end{array}$ 
is an equivalence relation.
\end{lemma}
\hfill $\Box$

Denote the equivalence = bordism class of $(M,g)$ by $[M,g]$.

\begin{lemma}
$[M \cup M', g \cup g'] = [M \# M', g \# g']$.
\end{lemma}
\hfill $\Box$

\me
\no
{\bf Remark.}
$(M \# M', g \# g')$ is metrically not uniquely defined but its
bordism class is.
\hfill $\Box$

\begin{lemma}
Set $[M,g] + [M',g'] = [M \cup M', g \cup g'] = [M \# M', g \# g']$.
Then $+$ is well defined and the set of all $[M^n,g]$ becomes an
abelian semigroup.
\end{lemma}
\hfill $\Box$

Denote by $\Omega^{nc}_n$ the corresponding Grothendieck group which
is the bordism group of all oriented open complete Riemannian
manifolds. Here $0$ is generated by the diagonal $\Delta$ and
$-[[M,g],[M',g']] = [[M',g'],[M,g]]$.

There is no reasonable approach for a calculation of $\Omega^{nc}_n$
known to us. $\Omega^{nc}_n$ is much to large. The situation rapidly
changes if we consider several refinements of the notion of bordism,
combining this with a component in ${\cal M}^n(mf)$ and having
additional conditions in mind, e.g. geometric conditions as 
nonexpanding ends or spectral conditions. Moreover, we would be
interested to have a geometric  realization of $0$ and $-[M,g]$.

First we consider bordism with compact support. Here we require
as above 1), 2) and aditionally

3) (cs). There exists a compact submanifold $C^{n+1} \subset B^{n+1}$
such that $B \setminus \intt C$ is a product bordism, i. e. 
$(B \setminus \intt C, g_B|_{B \setminus \intt C}) = ((M \setminus \intt C)
 \times [0,1],g|_{M \setminus \intt C} + dt^2)$.
3) (cs) implies 3. (after a compact change of the metric).

Write 
$\begin{array}{c} {} \\[-1ex] \sim \\[-1ex] b,cs \end{array}$ 
for the corresponding bordism. The corresponding bordism group will
be denoted by $\Omega^{nc}_n(cs)$. At the first glance, the calculation
of $\Omega^{nc}_n(cs)$ or at last the characterization of the bordism
classes seems to be very difficult. But this is not the case as we
indicate now. We connect $\Omega^{nc}_n(cs)$ with the components of
${}^b \comp{L,iso,rel}(\cdot) \subset {}^b {\cal M}^n_{L,iso,rel}(mf)$.

\me
\no
{\bf Remark.}
If $(M_1,g_1), (M_2,g_2) \in {}^b \comp{L,iso,rel}(M,g)$ then in 
general $(M_1,g_1) \# (M_2,g_2) \not\in {}^b \comp{L,iso,rel}(M,g)$.
\hfill $\Box$

Consider ${}^b \comp{L,iso,rel}(M,g)$, 
$\{ [M',g']_{cs} | (M',g') \in {}^b \comp{L,iso,rel}(M,g) \}$ 
and the subgroup
$\Omega^{nc}_n (cs, {}^b \comp{L,iso,rel}(M,g)) \subset \Omega^{nc}_n (cs)$
generated by
$\{ [M',g']_{cs} | (M',g') \in {}^b \comp{L,iso,rel}(M,g) \}$.
We know $\Omega^{nc}_n (cs)$ completely if we know all
$\Omega^{nc}_n (cs, {}^b \comp{L,iso,rel}(M,g))$
and we know 
$\Omega^{nc}_n (cs$, ${}^b \comp{L,iso,rel}(M,g))$
completely if we know
$\{ [M',g']_{cs} | (M',g') \in {}^b \comp{L,iso,rel}(M,g) \}$.
But the elements of 
$\{ [M',g']_{cs} | (M',g') \in {}^b \comp{L,iso,rel}(M,g) \}$
can be completely characterized by characteristic numbers which
we define now.

Fix $(M,g) \in {}^b \comp{L,iso,rel}(M,g)$, $M$ oriented.
Let $(M_1,g_1) \in {}^b \comp{L,iso,rel}(M,g)$ and
$\Phi : M \setminus K \lra M_1 \setminus K_1$ be an orientation
preserving isometry. Define Stiefel Whitney numbers of the
pair $(M_1,M)$ by
\[ w^{r_1}_1 \dots w^{r_n}_n (M_1,M) 
   := \sk{w^{r_1}_1 \dots w^{r_n}_n}{[K_1]} 
   + \sk{w^{r_1}_1 \dots w^{r_n}_n}{[K]} . \]
Similarly we define for $(M_1,M)$ and $n=4k$ Pontrjagin numbers
\[ p^{r_1}_1 \dots p^{r_k}_k (M_1,M) 
   := \int\limits_{K_1} p^{r_1}_1 \dots p^{r_k}_k (M_1)
   - \int\limits_K p^{r_1}_1 \dots p^{r_k}_k (M_1) \]
and the signature
\[ \sigma(M_1,M) := \sigma(K_1) + \sigma(-K) . \]

\begin{lemma}
$w^{r_1}_1 \dots w^{r_n}_n (M_1,M)$, $p^{r_1}_1 \dots p^{r_k}_k (M_1,M)$
and $\sigma(M_1,M)$ are well defined and
\bea
 w^{r_1}_1 \dots w^{r_n}_n (M_1,M) &=& 
 \sk{w^{r_1}_1 \dots w^{r_n}_n (K_1 \cup K)}{[K_1 \cup K]}, \non \\
 p^{r_1}_1 \dots p^{r_k}_k (M_1,M) &=&
 \sk{p^{r_1}_1 \dots p^{r_k}_k (K_1 \cup -K)}{[K_1 \cup -K]}, \non \\
 \sigma(M_1,M) &=& \sigma(K_1 \cup -K). \non
\eea
Here $K_1 \cup -K$ means $K_1 \begin{array}{c} {} \\[-1ex] \cup \\[-1ex] 
\Phi |_{\partial K} \end{array} -K$.
\end{lemma}

We refer to [9] for the very simple proof.
\hfill $\Box$

A complete characterization of bordism classes is now given by

\begin{theorem}
Fix $(M_1,g_1), (M_2,g_2) \in {}^b \comp{L,iso,rel}(M,g)$.
Then $(M_1,g_1) \begin{array}{c} {} \\[-1ex] \sim \\[-1ex] b,cs 
\end{array} (M_2,g_2)$ if and only if all characteristic numbers
of $(M_1,M)$ coincide with the corresponding characteristic numbers
of $(M_2,M)$.
\end{theorem}

We refer to [9] for the proof.
\hfill $\Box$

\begin{corollary}
The characterization of all elements of $\Omega^{nc}_n(cs)$
reduces to ''counting'' the (generalized) components of
 ${}^b {\cal M}^n_{L,iso,rel}(mf)$.
\end{corollary}
\hfill $\Box$

\me
\no
{\bf Remark.} 
If we restrict to closed oriented $n$--manifolds then 
${}^b \comp{L,iso,rel}(S^n)$ contains all closed oriented
$n$--manifolds, independent of any choice of the Riemannian
metric, and characteristic numbers defined above of
$(M^n,S^n)$ coincide with the characteristic numbers of
$M^n$. This follows from the definition above and the fact
that they vanish for $S^n$ or by cutting isometric collared
small discs from $M^n$, $S^n$ and gluing.
\hfill $\Box$

\me
\no
{\bf Examples.} 

1) Consider $M' = (P^{2k} \C \setminus$ int (small disc)
$\cup$ metric cylinder $S^{2k-1} \times [0, \infty [$, 
corresponding metric) and $M = (2k$--disc $\cup$ metric
cylinder $S^{2k-1} \times [0, \infty [$, corresponding
metric). Then $(M',g') \in {}^b \comp{L,iso,rel}(M,g)$
but $(M',g')$ is not cs--bordant to $(M,g)$ since
$\sigma(M',M)=1$, $\sigma(M',M')=0$.

2) For any $(M,g)$ there is a map $\Phi_n : \Omega_n \lra \{ 
[M',g']_{cs} | (M',g') \in {}^b \comp{L,iso,rel}(M,g) \}$
given by $[N] \lra [M \# N, g']_{cs}$ independent of
the metric on $N$.
\hfill $\Box$

There are many other types of bordism which are discussed
in [9]. We present still a type of bordism where $0$ and
$-[M,g]$ are geometrically realized. For this sake, we restrict
to metrics and bordisms of bounded geometry.

Let $(M^n,g)$, $(M'^n, g')$ be open, oriented, satisfying $(I)$
and $(B_k)$. We say $(M^n,g)$ is $(I)$, $(B_k)$--bordant
to $(M'^n, g')$ if there exists an oriented manifold $(B^{n+1},g_B$
with boundary $\partial B$ such that the following holds.

1) $(\partial B, g_B|_{\partial B}) = (M,g) \cup (_M',g')$.

2) There exists $\delta > 0$ such that $\exp: U_\delta(0_\nu)
\lra U_\delta(\partial B)$ is a $(k+2)$--bibounded diffeomorphism,
i. e. there exists a ''uniform collar'' of $\partial B$. Here
$0_\nu$ denotes the zero section of the inner normal bundle
$\nu=\nu(\partial B)$ of $\partial B$ in $B$.

3) $g_B$ satisfies $(B_k)$ on $B$ and $(I)$ on $B \setminus
U_{\frac{\delta}{2}} (\partial B)$.

4) There exists $R>0$ such that $B \subseteq U_R(M)$ and
$B \subseteq U_R(M')$. We write then $(M,g) \begin{array}{c} {} \\[-1ex]
\sim \\[-1ex] b,bg \end{array}$ $(M',g')$.

\begin{lemma}
Assume $(M,g) \begin{array}{c} {} \\[-1ex] \sim \\[-1ex] b,bg 
\end{array} (M',g')$ via $(B,g_B)$. Then there exist $\delta_1>0$
and $\tilde{g}_B$ s. t. $(M,g) \begin{array}{c} {} \\[-1ex]
\sim \\[-1ex] b,bg \end{array} (M',g')$ via
$(B,\tilde{g}), \tilde{g}_B |_{U_{\delta_1}(\partial B)} \cong
g_{\partial B} +dt^2$.
\end{lemma}

We refer to [9] for the proof.
\hfill $\Box$

\begin{corollary}
Without loss of generality we can always assume that the collar of
$\partial B$ is a metric collar.
\end{corollary}
\hfill $\Box$
\begin{corollary}
$\begin{array}{c} {} \\[-1ex] \sim \\[-1ex] b,bg \end{array}$ is an
equivalence relation.
\end{corollary}
\hfill $\Box$

Denote by $[M,g]_{bg}$ the equivalence = bordism class of $(M,g)$.

\begin{lemma}
$[M \cup M',g \cup g']_{bg} =[M \# M',g \# g']_{bg}$. Set 
$[M,g]_{bg} + [M',g']_{bg} := [M \cup M',g \cup g']_{bg} = 
[M \# M',g \# g']_{bg}$. Then $+$ is well defined and 
$\{ [M,g]_{bg} | (M,g)$ open, oriented with $(I)$ and $(B_k) \}$
becomes an abelian semigroup.
\hfill $\Box$
\end{lemma}

\begin{lemma}
Let $(M^n,g)$ be an open manifold satisfying $(I)$ and let
$\epsilon$ be an end of $M$. Then there exists a geodesic ray
$c$ tending to $\infty$ in $\epsilon$ with a uniformly thick
tubular neighborhood.
\end{lemma}
\hfill $\Box$

We define an end $\epsilon$ is nonexpanding ($\epsilon$ is a 
n. e. end), if there exists an $R>0$ and a ray $c$ in $\epsilon$
so that $\epsilon \subseteq U_R(|c|)$, which means that all 
elements of a neighborhood basis of $\epsilon$ are contained in 
$U_R(|c|)$.

We now restrict to manifolds with finitely many n. e. ends only.
Define 
$chc(r) := (D^n \cup S^{n-1}_r \times [0,\infty[, g_{standard})$, 
where $g_{standard}|_{S^{n-1}_r \times [a,\infty[}=
g_{S^{n-1}_r} +dt^2$ and the
standard metrics near $\partial D^n \begin{array}{c} {} \\[-2ex] 
= \\[-1ex] identif. \end{array} S^{n-1}_r \times \{ 0 \}$
are smoothed out. Then $chc(r)$ has one end, nonexpanding, 
and satisfies $(I)$ and $(B_\infty)$.

\begin{lemma}
$chc(r_1) \begin{array}{c} {} \\[-1ex] \sim \\[-1ex] b,bg \end{array} 
chc(r_2)$, $\bigcup\limits^s_{\sigma=1} chc(r_\sigma) \begin{array}{c} 
{} \\[-1ex] \sim \\[-1ex] b,bg \end{array} chc(r)$.
\end{lemma}
\hfill $\Box$ 

Define now 
\bea
  -[M,g] &:=& [-M,g], \non \\
  0 &:=& [chc(1)]
\eea
and set $\Omega^{nc,fe,ne}_n(I,B_k) := \all [M^n,g]_{bg}$ with
finitely many ends $(fe)$, all of them nonexpanding n.e..

\begin{theorem}
$(\Omega^{nc,fe,ne}_n(I,B_k), +,-,0)$ is an abelian group.
\end{theorem}

We refer to [9] for the proof, which contains throughout some 
delicate geometric constructions.
\hfill $\Box$

Examples of manifolds of type f.e., n.e., $(I)$, $(B_k)$ are
given by warped product metrics at infinity with $C^{k+2}$--bounded
warping function and e. g. by ends which are an infinite connected
sum of a finite number of closed Riemannian manifolds.

We finish with these two (very handable) examples
$\Omega^{nc}_n(cs)$, $\Omega^{nc,fe,ne}_n(I,B_k)$ our short
review of bordism theory for open manifolds and refer to [9]
for an extensive representation.

\section{Invariants of open manifolds}

\setcounter{equation}{0}

Consider $(M^n,g)$ with $(I)$ and $(B_k)$, usually $k \ge \frac{n}{2} +1$.
$(M^n,g)$ is a proper metric (length) space and we have sequences of
inclusions
\bea
&& \mbox{coarse type}(M,g) \supset \comp{GH}(M,g) \\
&& \mbox{coarse type}(M,g) \supset \comp{L}(M,g) \supset \comp{L,h}(M,g) 
\supset \comp{L,top}(M,g) \\
&& \comp{L}(M,g) \supset ^{b,m}\comp{L}(M,g) \supset ^{b,m}\comp{diff}(M,g) \\
&& ^{b,m}\comp{L}(M,g) \supset ^{b,m}\comp{L,h}(M,g) \supset ^{b,m}\comp{diff}(M,g) \\
&& ^b\comp{L,h,rel}(M,g) \supset ^b\comp{L,iso,rel}(M,g) 
\eea
and others. The arising task is to define for any sequence of inclusions
invariants depending only on the component and becoming sharper and sharper
if we move from left to right. We cannot present here all by us and others
defined invariants but present a certain choice. Start with the coarse type.
Denote by $C^*(M,g)$ the $C^*$--algebra obtained as the closure in ${\cal B}$
of $L_2{M}$ of all locally compact, finite propagation operators and by
$HX^*(M,g)$ the coarse cohomology. Then we have from [22], [23]

\begin{theorem}
$HX^*(M,g)$ and the $K$--theory $K_*(C^*(M,g))$ are invariants of the coarse
type, hence invariants of all components right from the coarse type.
\end{theorem}
\hfill $\Box$

\me
\no
{\bf Remarks.}

\no
1) We tried to describe the coarse type as the component of some uniform
structure but there is still some gap in this approach.

\no
2) Recall that a rough map is a coarse map which is uniformly proper. The
rough type is defined like the coarse type, replacing only coarse maps
by rough maps.
\hfill $\Box$

Block and Weinberger defined in [3] an homology $H^{Uf}_* (M,g)$ which
we here call rough homology (since it is functional under rough maps).
We prove in [8] that the rough and coarse type coincide.

\begin{theorem}
The rough homology $H^{Uf}_* (M,g)$ is an invariant of the coarse type.
\end{theorem}
\hfill $\Box$

Define the (singular) uniformly locally finite homology $H^{Uff}_* (M,g)$
as follows. It is the homology of the complex $C^{Uff}_* (M,g)$.
$c = \sum a_\sigma \sigma$ is a chain of $C^{Uff}_q$ if there exists $K>0$
depending on $c$ so that $|a_\sigma| \le K$ and the number of simplices
$\sigma$ lying in a ball of given size is uniformly bounded. The boundary
is defined to be the linear extension of the singular boundary. Similarly
is the uniformly locally finite cohomolgy $H^*_{Uff} (M,g)$ defined (cf.
[1]).

\begin{theorem}
$H^{Uff}_* (M,g)$ and $H^*_{Uff} (M,g)$ are invariants of $\comp{L,h}(M,g)$
\end{theorem}
\hfill $\Box$

Consider the bounded de Rham complex of $(M,g)$,
\[ \dots \lra ^{b,1}\Omega^q 
   \begin{array}{c} d \\[-1ex] \lra \\[-1ex] {} \end{array} 
   {}^{b,1}\Omega^{q+1} \lra \dots \]
Its cohomology is called the bounded cohomology of $^bH^* (M,g)$.

\begin{theorem}
$^bH^* (M,g)$ is an invariant of $^{b,2}\comp{L,h}(M,g)$.
\hfill $\Box$
\end{theorem}

\me
\no
{\bf Remark.}
We remember that we required in 5.4 the boundedness of the maps and the
homotopies. This is essential.
\hfill $\Box$

The $^b$--case is the smooth $L_\infty$--case and it is quite natural to
consider the $L_p$--case, in particular the $L_2$--case. But here the point
is that an arbitrary bounded map $f : (M,g) \lra (M',g')$ does not induce
an $L_2$--bounded map of analytic $L_2$--cohomology. Hence it is far from
being clear whether analytic $L_2$--cohomology is an invariant of 
$^{b,k+1}\comp{L,h}(M,g)$. This question can be attacked by simplicial
$L_2$--cohomology and $L_2$--Hodge--de Rham theory.

We will discuss this shortly and recall some simple definitions and facts.

Let $K$ be a locally oriented simplicial complex and $\sigma^q \in K$.
We denote $I(\sigma^q)= \# \{ \tau^{q+1} \in K | \sigma < \tau^{q+1} \}$,
$I_q(K) := \sup\limits_{\sigma^q \in K} I(\sigma^q)$. The complex $K$ is
called uniformly locally finite (u. l. f.) in dimension $q$ if
$I_q(K) < \infty$. If this holds for all $q$ then we call $K$ (in any
dimension) u. l. f.. The latter is equivalent to $I_0(K) < \infty$.
We assume in the sequel $K$ u. l. f.. Let for $1 \le p < \infty$ be
$C^{q,p}(K) = \{ c = \sum\limits_{\sigma^q \in K} c_\sigma \cdot \sigma |
c_\sigma \in \R, \sum\limits_\sigma |c_\sigma|^p < \infty \}$ the Banach
space of all real $p$--summable $q$--cochains. Then the linear extension
of $d$, 
$d \sigma^q = \sum\limits_{\tau^{q+1}} [\tau^{q+1} : \sigma^q] \tau^{q+1}$
is a bounded linear operator $d:C^{q,p} \lra C^{q+1,p}$,
$d(\sum c_\sigma \sigma) = \sum\limits_{\tau^{q+1}} 
\left( \sum\limits_{\sigma^q < \tau^{q+1}} [\tau : \sigma] \right) \tau^{q+1}$, and we
obtain a Banach cochain complex $(C^{*,p}, d)$. Its cohomology
$H^{*,p}(K)$ is called simplicial $L_p$--cohomology, 
$H^{q,p}(K) = Z^{q,p}/B^{q,p} = ker (d: C^{q,p} \lra C^{q+1,p})/
im (d: C^{q-1,p} \lra C^{q,p})$. 
$\overline{H}^{q,p}(K) = Z^{q,p} / \overline{B^{q,p}}$ 
is called reduced simplicial $L_p$--cohomology. For $p=2$, $C^{q,2}(K)$
is a Hilbert space via 
$\langle c, c' \rangle = \sum\limits_{\sigma^q} c_\sigma \cdot c'_\sigma$,
and we obtain a Hilbert complex $C^{*,2}(K)$. We refer to 
[17] and [2] for many   
simple proofs and interesting geometric examples.

Let $(C_{*,p} = C^{*,p}, \partial)$, 
$\partial \sigma^q = \sum\limits_{\tau^{q-1}} [\sigma^q : \tau^{q-1}] \tau^{q-1}$,
$\partial (\sum c_\sigma \sigma) = \sum\limits_{\tau^{q-1}} 
( \sum\limits_{\sigma^q > \tau^{q-1}} [\sigma : \tau] c_\sigma) \tau^{q-1}$,
be the Banach complex of $p$--summable real chains and $H_{*,p}(K)$ or
$\overline{H}_{*,p}(K)$ the corresponding $L_p$--homology or reduced
$L_p$--homology, respectively.

\begin{lemma}
If $H^{q,p} (K) \neq \overline{H}^{q,p}(K)$ then there exists an infinite
number of independent cohomology classes in $H^{q,p}$ whose image in 
$\overline{H}^{q,p}$ equals to zero, i. e. 
$\dim{\ker{H^{q,p} \lra \overline{H}^{q,p}}} = \infty$.
The same holds for $L_p$--homology.
\end{lemma}
\hfill $\Box$

\begin{corollary}
If $\dim{H^{q,p}(K)} < \infty$ then 
$H^{q,p}(K) = \overline{H}^{q,p}(K)$.
The same holds for homology.
\end{corollary}
\hfill $\Box$

\me
\no
{\bf Remark.}
$H^{q,p}(K)$ is endowed with a canonical topology, the quotient topology. 
But if $B^{q,p}$ is not closed then points in $H^{q,p}(K)$ are not closed.
In particular, any point $0 \neq \overline{b}+B^{q,p} \in H^{q,p}$,
$\overline{b} \in \overline{B}^{q,p} \setminus B^{q,p}$, belongs to
the closure $\overline{0} \subset H^{q,p}$.
\hfill $\Box$

For $p=2$ $\partial_{q-1}: C^{q,2} \lra C^{q-1,2}$ is the
$\langle , \rangle$--adjoint of 
$d_{q-1} : C^{q-1,2} \lra C^{q,2}$. Then
$\Delta_q (K) := d_{q-1} \partial_{q-1}+\partial_q d_q$
is a well defined bounded operator 
$\Delta_q : C^{q,2}(K) \lra C^{q,2}(K)$
and ${\cal H}^q(K) \equiv {\cal H}^{q,2}(K) := ker \Delta_q(K)$
is the Hilbert subspace of harmonic $L_2$--cochains = harmonic 
$L_2$--chains.

\begin{lemma}

\no
{\bf a)} $c \in {\cal H}^{q,2}(K)$ if and only if $dc=\partial c=0$

\no
{\bf b)} There exists an orthogonal decomposition
\[ C^{q,2} = {\cal H} ^q(K) \oplus dC^{q-1,2} \oplus \partial C^{q+1,2}\] 

\no
{\bf c)} There are canonical topological isomorphisms 
\[ {\cal H}^q \cong \overline{H}^{q,2}(K) \cong \overline{H}_{q,2}(K) . \]
\end{lemma}
\hfill $\Box$

Denote by $\sigma(\Delta_q(K))$ the spectrum, by $\sigma_e$ the
essential spectrum.

\begin{lemma}
The following conditions are equivalent.

\no
{\bf a)} $im \partial_q$ and $im \partial_{q-1}$ are closed.

\no
{\bf b)} $im d_q$ and $im d_{q-1}$ are closed.

\no
{\bf c)} $im \Delta_q(K)$ is closed.

\no
{\bf d)} $0 \notin \sigma_e (\Delta_q(K) |_{(ker \Delta_q)^\perp})$.

\no
{\bf e)} $H^{q,2}(K) = \overline{H}^{q,2}(K)$ and
$H^{q+1,2}(K) = \overline{H}^{q+1,2}(K)$.

\no
{\bf f)} $H_{q,2}(K) = \overline{H}_{q,2}(K)$ and
$H_{q-1,2}(K) = \overline{H}_{q-1,2}(K)$.
\end{lemma}
\hfill $\Box$

\me
\no
{\bf Example.}
Let $K=S^1 \times \R$ with a translation invariant u. l. f. triangulation.
Then $\overline{H}_{1,2}(K)=(0)$, $\dim{H_{1,2}(K)} = \infty$, i. e.
condition d) is not fulfilled.
\hfill $\Box$

For later applications the following simple lemma is quite useful.

\begin{lemma}
Let $H$ be a Hilbert space and $\Phi : H \lra H^{q,2}(K)$ be a vector space
isomorphism and homeomorphism. Then 
$H^{q,2}(K) = \overline{H}^{q,2}(K)$. The same holds as homology version.
\end{lemma}
\hfill $\Box$

We discuss next the behaviour of functional cohomology under maps and
subdivision.

We write $x(y)=\langle x, y \rangle$ for the value of a $q$--cochain $x$
and a $q$--chain $y$. If $S$ is a set of oriented $q$--simplexes then $S$
can be considered as a $q$--chain $S$ defined by 
$S = \sum\limits_{\sigma^q \in S} \sigma$.
Then $S$ will be called a geometric $q$--chain. Let $K$ and $L$ be
complexes and $T : C^{*,p}(K) \lra C^{*,p}(L)$ be a linear map. $T$ is
called vicinal if there exists an $N$ such that for all $\sigma \in K$
\[ \# \{ \tau \in L | \langle T \sigma, \tau \rangle \neq 0 \} \le N \]
and for any $\tau \in L$
\[ \# \{ \sigma \in K | \langle T \sigma, \tau \rangle \neq 0 \} \le N. \]

\begin{lemma}
If $T$ is vicinal and $| \langle T \sigma , \tau \rangle | \le M$ for all
$\sigma \in K$, $\tau \in L$ then $T$ is bounded.
\hfill $\Box$
\end{lemma}

Let $v$ be a vertex of $K$,
\[ st(v) = \{ \sigma \in K | v \mbox{ is a vertex of } \sigma \} , \]
and for $\sigma \in K$
\[ N(\sigma) = \overline{ \bigcup\limits_{v \in \sigma} st(v)} , \]
for $L \subset K$
\[ N(L) = \overline{ \bigcup\limits_{\sigma \in L} N(\sigma)} \]
\[ N^{(m)}(\cdot) := N( N( \dots (N(\cdot)) \dots ) . \]

We recall some definitions from [2].

A linear map $T: C^{*,p}(K_1) \lra C^{*,p}(K_2)$ is local if it is vicinal and 
there exists a positive integer $n$ such that whenever there is a
simplicial bijection $\eta: N^n(\sigma) \lra N^n(\tau)$ with 
$\eta(\sigma)=\tau$ then there exists a simplicial bijection
$\varepsilon : Im (N^n(\sigma)) \lra Im (N^n(\tau))$ such that
$T(\eta_*(\sigma)) = \varepsilon_*(T(\sigma))$.
Here $Im (N^n(\sigma))$ denotes the subcomplex supporting the chain
$T(N^n(\sigma))$. In other words, a local map is characterized by
the following three properties:

\no
{\bf 1)}
The image of a simplex $\sigma$ is a chain whose support has less than
$N$ simplexes, $N$ independent of $\sigma$.

\no
{\bf 2)}
The number of simplexes whose image contains a simplex $\tau$ is less
than $N$, $N$ independent of $\tau$.

\no
{\bf 3)}
The value of the map on a simplex $\sigma$ depends only on the
configuration arround $\sigma$.

These conditions are independent of $p$ and purely combinatorial in nature.

\me
\no
{\bf Example.}
$d$, $\partial$ and $\Delta$ are local operators.
\hfill $\Box$

For many applications one can weaken the third condition replacing the
existence of the map $\varepsilon$ by the weaker condition that there
exists a constant $c$ such that
\[ |T(\eta_*(\sigma))|_p < c \cdot |T(\sigma)|_p . \]
Such maps will be called nearly local.

\begin{lemma}
If $T$ is nearly local then it is bounded.
\end{lemma}

Now we want to prove the invariance of functional cohomology under 
subdivisions of finite degree and apply this for duality on open
manifolds. For this we need more general complexes than simplicial
ones. An absolute complex is a 3--couple $K=(K, <, \dim{})$, 
$\dim{}: K \lra \Z_+$, such that $<$ is a transitive relation, $\dim{}$ is
monotone w. r. t. $<$ and for every $x \in K$ there are only finitely
many $y \in K$ with $y<x$. If $\dim{x} = n$ then $x$ is called an 
$n$--cell. $\varepsilon : K \times K \lra \Z$ is called a boundary
function or incidence number if $\varepsilon(x,y) \neq 0$ implies
$x<y$ and $\dim{y} = 1 + \dim{x}$ and for $x,y \in K$, $\dim{y}=2+\dim{x}$
holds $\sum\limits_{z \in K} \varepsilon(x,z) \varepsilon(z,y)=0$.
We write $K_\varepsilon = (K,<,\dim{},\varepsilon)$. Then the definition
of u. l.f. $K_\varepsilon$ is quite clear and then the spaces
$C^{*,p}(K_\varepsilon,d), H^{*,p}(K_\varepsilon), 
\overline{H}^{*,p}(K_\varepsilon), C_{*,p}(K_\varepsilon), 
H_{*,p}(K_\varepsilon), \overline{H}_{*,p}(K_\varepsilon)$
are well defined. Simplicial, cell and CW--complexes produce their 
(after orientation) $K_\varepsilon$ . It is very easy and natural to
extend the notions local and nearly local to this general situation.

Let $K$ be an u. l. f. simplicial complex, $K'$ a subdivision. We say
that $K'$ has bounded degree of subdivision if for any $q$ there exists
a number $m_g$ s. t. any $\sigma^q \in K$ is subdivided into $\le m_q$
$q$--simplexes $\sigma'^q \in K'$. Then $K'$ is automatically u. l. f..
For $\sigma^q \in K$ let 
$B(\sigma^q)=\{ \sigma' \in K' | |\sigma'| \subset |\sigma^q| \}$.
Then ${\cal Z}(K') = \{ B(\sigma) | \sigma \in K \}$
is a cell decomposition of $K'$ (cf. [21], pp. 528--533 for definitions).
$B(\sigma)$ is a finite subcomplex and is called a block. There holds
\[ H_i(B(\sigma^q), \dot{B}(\sigma^q); \Z) \cong
   \left\{ \begin{array}{l} \Z, i=q \\ 0, \mbox{ in any other case.} 
   \end{array} \right. \]
Choose for any $B(\sigma^q)$ an integer orientation $b$, i. e. a generator
of $H_q(B, \dot{B}; \Z)$. Then we obtain oriented blocks (B, b). $b$ has
a representation $B=\sum\limits^\kappa_{i=0} \gamma_i \sigma'^q_i$, where
$\sigma'^q_i$ are the $q$--simplexes of $B(\sigma^q)$, $\gamma_i = \pm 1$,
$\kappa \le m_q$. Let $(B_1, b_1), \dots, (B_{q+1}, b_{q+1})$ be the
$(q-1)$--blocks of $\dot{B}$. Then $\partial b$ has a representation
$\partial b = \sum\limits^{q+1}_{\nu=1} \varepsilon_\nu b_\nu$, 
$\varepsilon_\nu = \pm 1$. $[(B,b) : (B_\nu,b_\nu)]' := \varepsilon_nu$
defines incidence numbers and a boundary operator.

Then it is possible by standard procedures to define a cell complex
$({\cal Z}(K'), [:]')$ such that the following holds

\begin{lemma}
After appropriate choice of orientation $b$, the map
$\sigma^q \lra B(\sigma^q)$ defines an isomorphism
\[ \zeta : (K, [:]) \lra  ({\cal Z} (K'), [:]') . \]
This holds even for arbitrary simplicial complexes and simplicial
subdivisions.
\end{lemma}
\hfill $\Box$

In the sequel we abbreviate ${\cal Z}(K') = ({\cal Z}(K'), [:]')$.

\begin{corollary}
The map $\zeta: \sigma^q \lra B(\sigma^q)$ induces canonical topological
isomorphisms $\zeta_{\#} : C_{*,p}(K) \lra C_{*,p}({\cal Z}(K'))$ and
$\zeta_* : H_{*,p}(K) \lra H_{*,p}({\cal Z}(K'))$. The same holds for
$C^{*,p}$ and $H^{*,p}$ and for the reduced case.
\end{corollary}
\hfill $\Box$

Let $(B(\sigma^q), b)$ be an elementary $q$--chain of of 
$C_{q,p}({\cal Z}(K'))$. Then we define
\[ \Theta_q (B(\sigma^q),b) := b = \sum\limits^\kappa_{i=1} \gamma_i
   \sigma'_i \in C_{q,p}(K'), \quad \kappa \le m_q . \]
Here $\sigma'^q$ are the $q$--simplexes of $B(\sigma^q)$ and 
$\gamma_i = \pm 1$. We define for 
$c=\sum\limits_{\sigma^q \in K} c_{B(\sigma^q)} (B(\sigma^q),b) \in 
C_{q,p}({\cal Z}(K'))$ the value $\Theta_qc$ by linear extension, 
$\Theta_qc=\sum\limits_{\sigma'^q} c_{\sigma'^q} \sigma'^q$,
where $c_{\sigma'^q}=\gamma_\nu \cdot c_{B(\sigma^q)}$,
$\sigma'^q \in B(\sigma^q)$.

\begin{lemma}
The map $\Theta=(\Theta_q)_q : C_{*,p}({\cal Z}(K')) \lra C_{*,p}(K')$
is an $L_p$--chain map. Any $\Theta_q : C_{q,p}({\cal Z}(K')) \lra C_{q,p}(K')$
is a monomorphism.
\hfill $\Box$
\end{lemma}

\begin{lemma}
$\Theta_q (C_{q,p}({\cal Z}(K')))$ is closed in $C_{q,p}(K')$.
\end{lemma}
\hfill $\Box$

\begin{corollary}
$\Theta_q : C_{q,p}({\cal Z}(K')) \lra im(\Theta_q)$ is a topological
isomorphism.
\end{corollary}
\hfill $\Box$

Define $\eta: K' \lra K$ as follows. For a vertex $v' \in K'$ denote by
$tr(e')$ the simplex of smallest dimension which contains $v'$. Choose
for $v'$ a vertex $v=\eta(v')$. $\eta$ is called a vertex translation
and defines by extension a simplicial map $\eta: K' \lra K$. Then the 
following lemma immediately follows from the corresponding lemma in
classical homology theory (cf. [21], 49.15, p. 537). We identify
$C_{q,p}(K')$ and $C_{q,p}({\cal Z}(K'))$ by means of $\zeta$ according
to 6.13.

\begin{lemma}
$(\eta_{\#q}) \circ \Theta_q : C_{q,p}(K) \lra C_{q,p}(K)$ equals
$id_{C_{q,p}(K)}$.
\end{lemma}
\hfill $\Box$

Now we can establish

\begin{theorem}
Let $1 \le p < \infty$ and $K$ be an u. l. f. simplicial complex and $K'$
a simplicial subdivision which is of bounded degree of subdivision in any
dimension. Then, after identification of $K$ with ${\cal Z}(K')$ by means
of $\zeta$, $\Theta : C_{*,p}(K) \lra C_{*,p}(K')$ induces topological
isomorphisms
\[ \Theta_* : H_{*,p}(K) \lra H_{*,p}(K') \]
and
\[ \Theta_* : \overline{H}_{*,p}(K) \lra \overline{H}_{*,p}(K') . \]
\end{theorem}
\hfill $\Box$

A first approach to homotopy invariance is given by

\begin{theorem}
Let $K$, $L$ be u. l. f. simplicial complexes, $1 \le p < \infty$
and $\varphi_{\#}: C_{*,p}(K) \lra C_{*,p}(L)$,
$\psi_{\#}: C_{*,p}(L) \lra C_{*,p}(K)$ chain maps satisfying
the following conditions.

(C1) \qquad $\varphi_{\#}$, $\psi_{\#}$ are local.

(C2) \qquad There are bounded chain homotopies
\[ \psi_{\#} \varphi_{\#} 
   \begin{array}{c} {} \\[-1ex] \sim \\[-1ex] D_1 \end{array} 
   id_{C_{*,p}(K)}, \qquad
   \varphi_{\#} \psi_{\#} 
   \begin{array}{c} {} \\[-1ex] \sim \\[-1ex] D_2 \end{array} 
   id_{C_{*,p}(L)} . \]
Then $\varphi_{\#}$ and $\psi_{\#}$ induce topological isomorphisms
\[ H_{*,p}(K) 
   \begin{array}{c} \psi_* \\[-1ex] \longleftarrow \\[-1ex] \cong \\[-1ex] 
   \lra \\[-1ex] \varphi_* \end{array} 
   H_{*,p}(L) \]
\[ \overline{H}_{*,p}(K) 
   \begin{array}{c} \psi_* \\[-1ex] \longleftarrow \\[-1ex] \cong \\[-1ex] 
   \lra \\[-1ex] \varphi_* \end{array} 
   \overline{H}_{*,p}(L) \]
\end{theorem}
\hfill $\Box$

Now we apply this to open manifolds $(M^n,g)$ satisfying $(I)$, $(B_k)$,
$k > \frac{n}{p} + 1$. Then it is a well known fact that such an $M$ 
admits a so called uniform triangulation, i. e. a triangulation 
$t: |K| \lra M$ with the following properties.

Let $\sigma^n$ be a curved $n$--simplex in $M^n$. We define the fullness
$\Theta (\sigma)$ by $\Theta(\sigma) = \vol (\sigma) / (\diam (\sigma))^n$.

{\bf a)} The exists a $\Theta_0 > 0$ such that for any curved simplex 
$\sigma^n$ the fullness satisfies the inequality $\Theta(\sigma) > \Theta_0$.

b) There exist constants $c_2>c_1>0$ such that for every $\sigma^n$ we
have
\[ c_2 \le \vol (\sigma) \le c_1 . \]

{\bf c)} There exists a constant $c>0$ such that for every vertex $v \in K$ the
barycentric coordinate function $\varphi_v : M \lra \R$ satisfies the
condition $|\nabla \varphi_v| \le c$.

If one assumes a) then b) is equivalent to the existence of bounds 
$d_1>d_2>0$ with $d_2 \le \diam (\sigma) \le d_1$ for all $\sigma \in K$.
a) and b) are equivalent to the boundedness of the volumes from below and
the diameters from above.

We call triangulations which satisfy the condition a) -- c) uniform. 
Uniform triangulations are u. l. f.. A connection between combinatorial 
and analytical theory is given by the following fundamental theorem of 
Goldstein / Kuzminov / Shvedov (cf. [18]).

\begin{theorem}
Let $t : |K| \lra M$ be a uniform triangulation and $1 \le p < \infty$.
Then there exists a canonical topological isomorphism $w_*$ (essentially
induced by the Whitney transformation), 
\bea
  w_* &:& H^{p,*}(K) \lra H^{p,*}(M,g), \non \\
  w_* &:& \overline{H}^{p,*}(K) \lra \overline{H}^{p,*}(M,g). \non
\eea
\end{theorem}
\hfill $\Box$

\me
\no
{\bf Remark.}
The $p=2$ version for the reduced case has been established by
Dodziuk (cf. [4]).
\hfill $\Box$

Consider now $(M^n_1,g_1)$, $(M^n_2,g_2)$ with $(I)$ and $(B_k)$ and a
uniformly proper map $f \in C^{\infty,k+1}(M_1,M_2)$.

\begin{proposition}
Let $t_i : |K_i| \lra M_i$, $i=1,2$, be uniform triangulations. Then there
exists a uniform triangulation $K'_1$ which has bounded degree of
subdivision and a simplicial approximation $f' : K^1 \lra K_2$ of $f$ which
is at $L_2$--chain level local.
\hfill $\Box$
\end{proposition}

\begin{samepage}
\begin{theorem}
Let $K^n$ be an oriented u. l. f. combinatorial homology $n$--manifold.
Then there exist isomorphisms
\bea
  D &:& H^{*,2}(K^n) \lra H_{n-*,2}(K^n) \non \\
  \overline{D} &:& \overline{H}^{*,2}(K^n) \lra \overline{H}_{n-*,2}(K^n) \non
\eea
\hfill $\Box$
\end{theorem}
\end{samepage}

Combining now 6.19 -- 6.22, we obtain

\begin{theorem}
Assume $(M^n_i,g_i)$ with $(I)$ and $(B_k)$, $k>\frac{n}{2}+1$. Let
$f \in C^{\infty,k+1}(M_1,M_2)$, $g \in C^{\infty,k+1}(M_2,M_1)$
uniformly proper homotopy equivalences, inverse to each other. Then 
$f$ and $g$ induce topological isomorphisms
\bea
  H^{2,*}(M_1,g_1) & \begin{array}{c} \Phi \\[-1ex] \lra \\[-2ex] 
  \longleftarrow \\[-1ex] \Psi \end{array} H^{2,*}(M_2,g_2) \non \\
  \overline{H}^{2,*}(M_1,g_1) & \begin{array}{c} \Phi \\[-1ex] \lra \\[-2ex] 
  \longleftarrow \\[-1ex] \Psi \end{array} \overline{H}^{2,*}(M_2,g_2) \non 
\eea
where $\Phi$, $\Psi$ are induced by simplicial approximations of 
$f$, $g$, $\Theta$, $\eta$ and $w$.
\end{theorem}
\hfill $\Box$

\begin{corollary}
Assume $(M^n_i,g_i)$, $f$, $g$ as above, $0 \le q < n$. Then
$\inf \sigma_e (\Delta_q (M_1) |_{(\ker \Delta_q(M_1))^\perp}) > 0$
if and only if
$\inf \sigma_e (\Delta_q (M_2) |_{(\ker \Delta_q(M_2))^\perp}) > 0$.
\end{corollary}

\me
\no
{\bf Proof.}
Apply the analytical version of 5.8. Then the spectral gap exists for
$\Delta_q(M_1)$ if and only if 
$H^{2,q}(M_1,g_1) = \overline{H}^{2,q}(M_1,g_1)$,
$H^{2,q+1}(M_1,g_1) = \overline{H}^{2,q+1}(M_1,g_1)$.
But this holds if and only if it holds for $(M_2,g_2)$.
\hfill $\Box$

We will not discuss here the general theory of characteristic numbers
for open manifolds, $L_2$--intersection theory and general duality. We
only mention that for manifolds of bounded geometry duality holds from
$H^*_{uff}$ to $H^{uff}_*$ and $L_2$-bounded duality in the $L_2$--case.
This yields obstructions for the corresponding components to contain a
manifold of bounded geometry.

We also defined homologies for the rel--components but refer to
[8], [9], [10].

\bi
\bi
J\"urgen Eichhhorn

Institut f\"ur Mathematik und Informatik

Friedrich--Ludwig--Jahn--Stra\ss e 15a

D-17487 Greifswald

Germany

eichhorn@mail.uni-greifswald.de

\end{document}